\documentclass[10pt]{amsart}
\usepackage{epsfig}
\usepackage{graphicx}

\usepackage{epsf}
\usepackage{epic,eepic}

\newtheorem{theorem}{Theorem}[section]
\newtheorem{lemma}[theorem]{Lemma}

\newtheorem{proposition}[theorem]{Proposition}

\theoremstyle{definition}

\theoremstyle{remark}
\newtheorem{remark}[theorem]{Remark}

\numberwithin{equation}{section}

\newcommand{\Er}{\mathbb{R}}

\newcommand{\En}{\mathbb{N}}
\newcommand{\dt}{\Delta t}

\title[Vibrations of a beam between stops \ \   ]{Vibrations of a beam between stops: Convergence of a fully discretized approximation}

\author{Yves Dumont}
 \address{IREMIA\\
  Universit\'e de La R\'eunion \\
15 avenue R. Cassin\\
 97715 Saint-Denis Messag. 9, France}
\email{Yves.Dumont@univ-reunion.fr}

\author{ Laetitia Paoli}

\address{Equipe d'Analyse Num\'erique de Saint-Etienne \\
Universit\'e Jean Monnet \\
23 Rue du Docteur Paul Michelon\\
42023 St-Etienne Cedex 2, France}
\email{laetitia.paoli@univ-st-etienne.fr}
\begin{document}

\begin{abstract}
 We consider the dynamics of an elastic beam which is clamped at its left end to a vibrating support and which can move freely at its right end between two rigid obstacles (the stops). We model the contact with Signorini's complementary conditions between the displacement and the shear stress. 
For this infinite dimensional contact problem, we propose a family of fully discretized approximations and their convergence is proved. Moreover some examples of implementation are presented.

\vspace{0.5cm}

{\bf Keywords}: Dynamics with impact -- Signorini's conditions -- Space and time discretization -- Convergence.

{\bf AMS}: 35L85, 65M12, 74H45
\end{abstract}

\maketitle

\section{Description of the problem} \label{intro}

We consider a beam which is clamped at its left end to a vibrating support and which can move freely between two rigid obstacles at its right end (see figure 1).
 
The longitudinal axis of the beam coincide with the interval $[0,L]$ and we denote by $\tilde u (x,t)$, $(x,t) \in (0,L) \times (0,T)$ the vertical displacement of a point $x$ belonging to this axis. We assume that the material is elastic and the motion is planar. We denote by $\tilde \sigma$  the shear stress given by
\begin{eqnarray*}
\tilde \sigma (x,t) = - k^2 \tilde u_{xxx},\quad  k^2 = \frac{EI }{ \rho S} 
\end{eqnarray*}
where $\rho$ and $E$ are the density and the Young's modulus of the material and $S$ and $I$ are respectively the surface and the inertial momentum of the section of the beam.
Then, under the assumption of small displacements, the motion is described by the following partial differential equation
\begin{eqnarray*}
\tilde u_{tt} - \tilde \sigma_x = \tilde f 
\end{eqnarray*}
where $\tilde f$ is the density of external forces.
\begin{figure}
\begin{center}
\begin{picture}(10,55)
\put(-50,0){\begin{picture}(1,1)
\put(0,-9){\line(0,1){48}}
\put(0,10){\framebox(80,10){beam}}
\put(45,10) {\vector(0,-1){12}}
\put(35,-2){$u$}
\put(75,0){\framebox(8,5)}
\put(75,25){\framebox(8,5)}
\put(83,-9){\line(0,1){48}}
\put(-28, 12){$\phi$}
\put(-17,4){\vector(0, 1){20}}
\multiput(0,0)(0,5){7}{\line(-2,-1){10}} 
\put(110,15){stops}
\put(103,27){\vector(-1,0){20}}
\put(103,3){\vector(-1,0){20}}
\put(103,3){\line(0,1){24}}
\end{picture}}
\end{picture}

\vspace{0.7 cm}
\caption{ The physical setting}
\label{fig:1}  
\end{center} 
\end{figure}
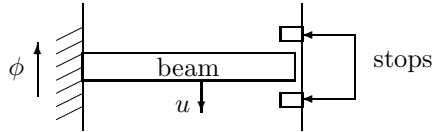

The beam is clamped at its left end so
\begin{eqnarray*}
\tilde u(0,t) = \phi (t), \quad \tilde u_x (0,t)=0 
\end{eqnarray*}
where $\phi$ describes the motion of the vibrating support. At its right end the beam can move freely between two obstacles, called "the stops", so we have
\begin{eqnarray*}
g_1 \le \tilde u (L,t) \le g_2, \quad \tilde u_{xx} (L, t) = 0  
\end{eqnarray*}
and we assume that $g_1 < 0 < g_2$. When the beam hits one of the two stops, the stress is in the opposite direction of  the displacement and we obtain the following Signorini's conditions
\begin{eqnarray*}
\left\{
\begin{array}{l}
 \tilde \sigma (L,t) \ge 0 \quad \hbox{\rm if $\tilde u(L,t)=g_1$,} \\ 
 \tilde \sigma (L,t) \le 0 \quad \hbox{\rm if $\tilde u(L,t)=g_2$,} \\  
 \tilde \sigma (L,t) =0  \quad \hbox{\rm if $g_1< \tilde u(L,t) < g_2$.} 
\end{array} \right.
\end{eqnarray*}

These relations can be rewritten as follows
\begin{eqnarray*}
 - \tilde \sigma (L,t) \in \partial \psi_{[g_1, g_2]} \bigl( \tilde u(L,t) \bigr)
\end{eqnarray*}
where $\psi_{[g_1, g_2]}$ is the indicator function of the interval $[g_1, g_2]$ and $\partial \psi_{[g_1, g_2]} $ is its subdifferential (\cite{rock}).

\smallskip

In order to deal with homogeneous boundary conditions at $x=0$, we consider a new unknown function $u$ defined by
\begin{eqnarray*}
u(x,t) = \tilde u (x,t) - h(x) \phi (t), 
\end{eqnarray*}
with
\begin{eqnarray*}
h(x) = 1 - 2 \left( \frac{x }{ L} \right)^2 + \frac{4 }{ 3} \left( \frac{x }{ L} \right)^3 - \frac{1 }{ 3} \left(\frac{x }{ L} \right)^4.
\end{eqnarray*}
The mechanical problem is now described by the system
\begin{eqnarray*}
\left\{ 
\begin{array}{l}
 u_{tt} + k^2 u_{xxxx} = f \quad \hbox{\rm in $(0,L) \times (0,T)$ } \\
 u(0,\cdot )=u_x (0,\cdot )= u_{xx} (L,\cdot )=  0 \quad \hbox{\rm in $(0,T)$} \\
 u(L,\cdot ) \in [g_1, g_2], \  u_{xxx} (L,\cdot ) \in \partial \psi_{[g_1, g_2]} \bigl(  u(L,\cdot ) \bigr) \quad \hbox{\rm in $(0,T)$}  
\end{array}\right. 
\end{eqnarray*}
with $f(x,t)=\tilde f(x,t) - h(x) \phi''(t) - k^2 h^{(4)} (x)  \phi(t) $ for all $(x,t) \in (0,L) \times (0, T)$. We complete the model with the initial conditions 
\begin{eqnarray*}
u(\cdot,0)=u_0, \quad u_t (\cdot ,0)=v_0 \quad \hbox{\rm in $(0, L)$.}
\end{eqnarray*}

\smallskip

As usual in mechanical problems with unilateral constraints we cannot expect classical solutions since the velocities may be discontinuous. So we look for weak solutions. For this purpose we introduce the following functional spaces 
\begin{eqnarray*}
H=L^2(0,L), \quad V= \bigl\{ w\in H^2(0,L) ; w(0)=w_x(0)=0 \bigr\}, 
\end{eqnarray*}
\begin{eqnarray*}
{\mathcal H} = \big\{ w \in L^2 (0,T; V); w_t \in L^2(0,T; H) \bigr\},
\end{eqnarray*}
and the convex set
\begin{eqnarray*}
 K=\bigl\{ w \in V; g_1 \le w(L) \le g_2 \bigr\}  .
\end{eqnarray*}
We denote by $(.,.)$ and $|.|$ the canonical scalar product and norm of $H$.
Let $a$ be the following bilinear form
\begin{eqnarray*}
a(u,v)= \int_0^L k^2 u_{xx} v_{xx} \, dx \quad \forall (u,v) \in V^2.
\end{eqnarray*}
We may observe that $a$ defines a scalar product on $V$ and the associated norm, denoted $\|.\|_V$, is equivalent to the canonical norm of $H^2(0,L)$ on $V$.
The weak formulation of the problem is then given by the following variational inequality
\begin{eqnarray*}  
(P) \left\{
\begin{array}{l}
\displaystyle  - \int_0^T \bigl(u_t (\cdot,t), w_t (\cdot, t) -u_t (\cdot, t) \bigr) \, dt +  \int_0^T a\bigl(u (\cdot,t), w (\cdot, t)-u (\cdot, t) \bigr) \\
\displaystyle   \ge \bigl(v_0, w(\cdot, 0) - u_0 \bigr) + \int_0^T \bigl(f(\cdot,t), w (\cdot,t) - u (\cdot,t) \bigr) \, dt \\
 \forall w \in {\mathcal H} \cap L^2(0,T; K) \ \hbox{\rm such that $w( \cdot, T) = u (\cdot, T)$.}
\end{array}
\right.
\end{eqnarray*}
  
For this problem an existence result has been obtained by K.Kuttler and M.Shillor by using a penalty method.

\smallskip

\begin{theorem} \label{theorem1}
 (\cite{Ku-Shi}) Assume that $f \in L^2(0,T; H)$, $u_0 \in K$, $v_0 \in H$. Then there exists $u \in {\mathcal H} \cap L^2(0,T; K)$ such that problem (P) is satisfied and $u(\cdot, 0)=u_0$.
\end{theorem}

\smallskip

 It should be noted that, as far as we know,  uniqueness remains an open question.

\smallskip

For the computation of approximate solutions, the penalty method which is introduced as a theoretical tool to obtain existence in \cite{Ku-Shi} could appear as an interesting technique: the Signorini's conditions are replaced by a normal compliance law  
\begin{eqnarray*}
\sigma (L,t) = - \frac{1}{\varepsilon} \bigl[ \max\bigl( u(L,t )-g_2, 0 \bigr) -  \max\bigl(g_1 -  u(L,t ) , 0 \bigr) \bigr], \quad \frac{1}{\varepsilon}>>1 
\end{eqnarray*}
which leads to a system of partial differential equations depending on the penalty parameter $\varepsilon$. From the mechanical point of view $1/\varepsilon$. can be interpreted as the stiffness of the stops which are not assumed to be perfectly rigid anymore. From a numerical point of view, for large values of $1/\varepsilon$, we have to solve a stiff problem which is expensive (\cite{Du}). Moreover the dynamics of the system may be complex (see \cite{ShMo} for a periodic forcing) and the approximate motion could be quite sensitive with respect to the value of $1/\varepsilon$. (see \cite{Barc} for an example in the case of a simplified model of vibrations, see also \cite{Bro}). 

In order to avoid these difficulties, we propose to deal directly with the unilateral boundary condition by solving a complete discretization, in both time and space, of the variational inequality (P). 

\smallskip

From now on we will consider the more general case of a convex set $K$ given by
\begin{eqnarray*}
K= \bigl\{ w \in V; g_1(x) \le w(x) \le g_2(x) \quad \forall x \in [0,L] \bigr\}
\end{eqnarray*}
where $g_1$, $g_2$ are two mappings from $[0,L]$ to ${\overline \Er}$ such that there exists $g>0$ such that
\begin{eqnarray*}
g_1(x) \le -g < 0< g \le g_2(x) \quad \forall x \in [0,L].
\end{eqnarray*}
Problem (P)  then  describes the vibrations of an elastic beam between two longitudinal rigid obstacles. 

\smallskip

The paper is organized as follows: in the next section we introduce the fully discretized approximation of the problem, then in section \ref{sec:3}, we prove its stability and convergence and finally, in section \ref{sec:4}, we present some examples of implementation. 

Let us observe that the convergence result yields also an existence result for the more general case that we consider here.

\section{Discretization} \label{sec:2}

Let us  assume now that the assumptions of theorem \ref{theorem1} hold i.e $f \in L^2(0,T; H)$, $u_0 \in K$ and $v_0 \in H$. 
\smallskip

For all $h \in \Er^*_+$ we consider a finite dimensional subspace $V_h$  of $V$  such that,
for all $v \in V$, there exists a sequence $(v_h)_{h>0}$ such that
\begin{eqnarray*}
\|v_h - v\|_V \to_{h \to 0} 0, \quad v_h \in V_h \quad \forall h >0,
\end{eqnarray*}
and we denote by $Q_h$ the projection onto $V_h$ respectively to the scalar product defined by $a$ on  $V$. The compact embedding of $V$ into $H^1(0,L)$ implies that 
  there exists a sequence $(\gamma_h)_{h >0}$ such that
\begin{eqnarray*}
\forall w \in V \quad \bigl\| Q_h (w) - w \bigr\|_{H^1(0,L)} \le \gamma_h \| w\|_V , \quad \lim_{h \to 0} \gamma_h =0.
\end{eqnarray*}
For all $h> 0$ we define $K_h = K \cap V_h$.

\bigskip

Let $N \in \En^*$ and $\Delta t = T/N$. We propose the following family of   discretizations of problem (P): For all $n \in \{ 1, \dots, N-1 \}$, find $u_h^{n+1} \in K_h$ such that
\begin{eqnarray*}
(P_{h \beta}^{n+1}) 
\left\{
\begin{array}{l}
\left(\displaystyle  \frac{u_h^{n+1} -2 u_h^{n} + u_h^{n-1} }{ \Delta t^2}, w_h - u_h^{n+1} \right) \\
\displaystyle  + a \left(  \beta u_h^{n+1} + (1-2 \beta)  u_h^{n} + \beta u_h^{n-1} , w_h - u_h^{n+1} \right) \\
\ge \left( \beta f^{n+1} + (1-2 \beta) f^{n} +  \beta f^{n-1}  , w_h - u_h^{n+1} \right) \quad \forall w_h \in K_h  
\end{array}
\right.  
\end{eqnarray*}
with 
\begin{eqnarray*}
f^{n} = \frac{1 }{ \dt} \int_{n \dt}^{(n+1) \dt} f( \cdot, s) \, ds   
\end{eqnarray*}
where $\beta$ is a parameter belonging to $[0, 1/2]$. We choose $u_h^0$ and $u_h^1$ in $K_h$ such that $\bigl( \| u_h^1 \|_V\bigr)_{h>0}$ remains bounded and 
\begin{equation} \label{eq0}
 \lim_{h \to 0, \Delta t \to 0} \| u_h^0 -u_0 \|_V + \left| \frac{u_h^1 - u_h^0 }{\dt} - v_0  \right| =0. 
\end{equation}

\bigskip

We may observe that problem $(P_{h \beta}^{n+1})$ ($\beta \in [0,1/2]$) can be rewritten as 
\begin{eqnarray*}
\left\{
\begin{array}{l}
 \hbox{\rm Find $u_h^{n+1} \in K_h$ such that} \\
 a_{n \beta } \bigl(u_h^{n+1}, w_h - u_h^{n+1} \bigr) \ge L_{n \beta} \bigl(w_h - u_h^{n+1} \bigr) \quad \forall w_h \in K_h
 \end{array} \right.
 \end{eqnarray*}
with
\begin{eqnarray*}
\begin{array}{l} 
 a_{n \beta} (u_h, v_h) = (u_h, v_h) + \dt^2 \beta  a (u_h, v_h) ,  \\
 L_{n \beta}(v_h) =  \dt^2 \bigl( \beta f^{n+1} + (1-2 \beta) f^{n} + \beta f^{n-1}, v_h \bigr)  + ( 2 u_h^{n} - u_h^{n-1}, v_h )  \\
 \qquad \quad -  \dt^2  a \bigl( (1-2 \beta) u_h^{n} +  \beta u_h^{n-1}, v_h \bigr)    
 \end{array}
 \end{eqnarray*}
for all $(u_h, v_h) \in V_h^2$.

By an immediate induction on $n$, we obtain that $L_{n \beta}$ is linear and continuous on $V_h$, and it is obvious that $a_{n \beta }$ is bilinear, symmetric, continuous and coercive on $V_h$. Thus the existence and uniqueness of $u_h^{n+1}$ follows.

\bigskip

\begin{remark} This family of discretizations is inspired by Newmark's algorithms of parameters $\gamma=1/2$, $\beta \in [0, 1/2]$ (\cite{Hu}). Indeed, if $K=V$ (i.e $g_1(x) = - \infty$, $g_2(x) = + \infty$ for all $x \in [0,L]$), the mechanical problem is described by the following system of partial differential equations
\begin{eqnarray*}
\left\{
\begin{array}{l}
 u_{tt} + k^2 u_{xxxx} = f \quad \hbox{\rm in $(0,L) \times (0,T)$} \\
 u (0, \cdot) = u_x (0, \cdot) =  u_{xx}(L, \cdot) = u_{xxx} (L, \cdot) =0  \quad \hbox{\rm in $(0,T)$} 
 \end{array} \right.
 \end{eqnarray*}
with the initial data $(u_0, v_0)$, and 
 $(P_{h \beta}^{n+1})$ reduces to 
\begin{eqnarray*}
\left\{
\begin{array}{l}
 \hbox{\rm find $u_h^{n+1} \in K_h$ such that } \\
 \left(\displaystyle  \frac{u_h^{n+1} -2 u_h^{n} + u_h^{n-1} }{ \Delta t^2}, w_h  \right) + a \left(  \beta u_h^{n+1} + (1-2 \beta)  u_h^{n} + \beta u_h^{n-1} , w_h  \right) \\
= \left( \beta f^{n+1} + (1-2 \beta) f^{n} +  \beta f^{n-1}  , w_h  \right) \quad \forall w_h \in V_h 
\end{array} \right. 
\end{eqnarray*}
which is simply a Newmark's scheme of parameters $\gamma = 1/2$, $\beta$ for the previous system.
\end{remark}

\section{Convergence} \label{sec:3}

Since the proposed  discretizations are inspired by  Newmark's methods which stability depends on the value of $\beta$, we may expect the same kind of result for $(P_{h \beta}^{n+1}) $.  More precisely, for $\beta \in [0, 1/2)$ we  obtain the following conditional  stability property:

\bigskip

\begin{proposition} \label{stab1}  Let $\beta \in [0, 1/2)$,  $h>0$  and $\kappa_h$ be defined by
\begin{eqnarray*}
\kappa_h = \sup_{u_h \in V_h\setminus \{0\} } \frac{a(u_h, u_h)}{ |u_h|^2 }.
\end{eqnarray*}
Let $\alpha \in (0,1)$ and $N_h \in \En^*$   be such that
\begin{equation} \label{eq1}
\frac{T }{N_h} < \min \left( 2 \sqrt{ \frac{1 - \alpha} { \kappa_h (1-2 \beta) } } , \alpha \right). 
\end{equation}
Then there exists a constant depending only on the data, $C( f, u_0, v_0)$, such that for all $h>0$ and for all $N \ge N_h$ (i.e $\displaystyle \dt < \dt_h= \frac{T }{ N_h}$)
\begin{eqnarray*}
\alpha \left| \frac{u_h^{n+1} - u_h^{n} }{ \dt } \right|^2 + \beta a \bigl( u_h^{n}, u_h^{n} \bigr) 
 + \beta a \bigl( u_h^{n+1}, u_h^{n+1} \bigr)\le C (f, u_0, v_0)
\end{eqnarray*}
for all $n \in \{ 1, \dots, N-1 \}$, where  $(u_h^{n+1})_{1 \le n \le N-1}$ are  the solutions of problems $(P_{h \beta}^{n+1})_{1 \le n \le N-1}$. 
\end{proposition}

\begin{remark} An estimate of $\kappa_h$ in the case of a P3 finite element space discretization is given in the Appendix.
\end{remark}

\begin{proof} Let $n \in \{ 1, \dots, N-1 \}$ and  choose $w_h = u_h^{n-1}$ as a test-function in $(P_{h \beta}^{n+1})$. We get
\begin{eqnarray*}
\left\{
\begin{array}{l}
\left(\displaystyle  \frac{u_h^{n+1} -2 u_h^{n} + u_h^{n-1} }{ \dt^2}, u_h^{n-1} - u_h^{n+1} \right) \\
+ a \left( \beta u_h^{n+1} + (1-2 \beta)  u_h^{n} + \beta u_h^{n-1} , u_h^{n-1} - u_h^{n+1} \right)
 \ge \bigl( g^{n} , u_h^{n-1} - u_h^{n+1} \bigr)
\end{array} \right.
\end{eqnarray*}
where 
\begin{eqnarray*}
g^{n} = \beta f^{n+1} + (1-2 \beta) f^{n} + \beta  f^{n-1}.
\end{eqnarray*}
The two first terms can be rewritten as follows:
\begin{eqnarray*}
\left( \frac{ u_h^{n+1} -2 u_h^{n} + u_h^{n-1} }{ \dt^2 }, u_h^{n-1} - u_h^{n+1} \right) 
= \left| \frac{u_h^{n-1} - u_h^{n}  }{ \dt } \right|^2 - \left| \frac{u_h^{n+1} - u_h^{n}  }{ \dt}  \right|^2 ,
\end{eqnarray*}
and
\begin{eqnarray*}
\begin{array}{l}
 a \left( \beta  u_h^{n+1} + (1 -2 \beta) u_h^{n} + \beta u_h^{n-1} , u_h^{n-1} - u_h^{n+1} \right) = (1-2 \beta)  a \bigl( u_h^{n-1}, u_h^{n} \bigr) \\
 + \beta a \bigl( u_h^{n-1}, u_h^{n-1} \bigr) 
 - (1 -2 \beta)  a \bigl( u_h^{n}, u_h^{n+1} \bigr)   - \beta a \bigl( u_h^{n+1}, u_h^{n+1} \bigr) .
 \end{array} 
 \end{eqnarray*}
Hence, for all $n \in \{ 1, \dots, N-1 \}$, we have
\begin{eqnarray*}
\begin{array}{l}
 \left|\displaystyle  \frac{u_h^{n+1} - u_h^{n} }{ \dt} \right|^2 + (1-2 \beta) a \bigl( u_h^{n}, u_h^{n+1} \bigr) + \beta a \bigl( u_h^{n+1}, u_h^{n+1} \bigr) \\
  \le \displaystyle  \left| \frac{u_h^{n} - u_h^{n-1} }{ \dt} \right|^2 + (1-2 \beta)   a \bigl( u_h^{n-1}, u_h^{n} \bigr) + \beta a \bigl( u_h^{n-1}, u_h^{n-1} \bigr)
   + \bigl( g^{n}, u_h^{n+1} - u_h^{n-1} \bigr)
\end{array}
\end{eqnarray*}
and with a discrete integration
\begin{eqnarray*}
\begin{array}{l}
 \left|\displaystyle  \frac{u_h^{n+1} - u_h^{n} }{ \dt} \right|^2 + (1-2\beta)  a \bigl( u_h^{n}, u_h^{n+1} \bigr) + \beta a \bigl( u_h^{n+1}, u_h^{n+1} \bigr) + \beta a \bigl( u_h^{n}, u_h^{n} \bigr) \\
\le \displaystyle \left| \frac{u_h^{1} - u_h^{0} }{ \dt} \right|^2 + (1-2 \beta) a \bigl( u_h^{0}, u_h^{1} \bigr) + \beta a \bigl( u_h^{1}, u_h^{1} \bigr) + \beta a \bigl( u_h^{0}, u_h^{0} \bigr)\\
\displaystyle  + \sum_{p=1}^n \bigl( g^p,  u_h^{p+1}- u_h^{p-1} \bigr).
\end{array}
\end{eqnarray*}
Using the same techniques as in \cite{MS-Ber}, we define
\begin{eqnarray*}
R(u_h, v_h) = (1-2 \beta) a(u_h, v_h) + \left| \frac{u_h - v_h }{ \dt} \right|^2  \quad \forall (u_h, v_h) \in V_h^2.
\end{eqnarray*}
We observe that
\begin{eqnarray*}
\begin{array}{l}
 R(u_h, v_h)  = \displaystyle  \frac{1 -2 \beta }{ 4} a (u_h + v_h, u_h + v_h) - \frac{1 - 2 \beta  }{ 4} a (u_h - v_h, u_h - v_h) 
\displaystyle  + \left| \frac{u_h - v_h }{ \dt} \right|^2 \\
\displaystyle  \ge \frac{1 -2 \beta }{ 4} a (u_h + v_h, u_h + v_h) 
\displaystyle  + \left( 1 - \kappa_h \dt^2 \frac{ 1-2\beta }{ 4} \right) \left| \frac{u_h - v_h }{ \dt} \right|^2 ,
 \end{array}
 \end{eqnarray*}
and with assumption (\ref{eq1}), we infer that
\begin{eqnarray*}
R(u_h, v_h) \ge \displaystyle  \frac{1 -2 \beta  }{ 4} a (u_h + v_h, u_h + v_h)  + \alpha \left| \frac{u_h - v_h }{ \dt} \right|^2 \quad  \forall (u_h, v_h) \in V_h^2 .
\end{eqnarray*}
It follows that
\begin{equation} \label{eq2}
\begin{array}{l}
\displaystyle   \alpha   \left| \frac{u_h^{n+1} - u_h^{n} }{ \dt} \right|^2 
+ \frac{1 -2 \beta }{ 4} a \bigl( u_h^{n+1} + u_h^{n} , u_h^{n+1} + u_h^{n} \bigr)
 + \beta a \bigl( u_h^{n+1}, u_h^{n+1} \bigr) + \beta a \bigl( u_h^{n}, u_h^{n} \bigr) \\
\displaystyle   \le R \bigl( u_h^0, u_h^1 \bigr) + \beta a \bigl( u_h^{1}, u_h^{1} \bigr)
 \displaystyle + \beta a \bigl( u_h^{0}, u_h^{0} \bigr) 
 + \sum_{p=1}^n \bigl| g^p \bigr|^2 \dt  + \sum_{p=0}^{n} \left| \frac{u_h^{p+1} - u_h^{p} }{ \dt} \right|^2 \dt .
\end{array} 
\end{equation}

Since $\alpha - \dt \ge \alpha - \dt_h >0$, Gr\"onwall's lemma implies that
\begin{eqnarray*}
\sum_{p=0}^n \left| \frac{u_h^{p+1} - u_h^{p} }{ \dt} \right|^2 \le \left| \frac{u_h^{1} - u_h^{0} }{ \dt} \right|^2 \exp \left( \frac{ n  \dt }{ \alpha - \dt} \right) 
+ \sum_{p=1}^n k_p \exp \left( \frac{ (n-p) \dt }{ \alpha - \dt }  \right)
\end{eqnarray*}
with
\begin{eqnarray*}
k_p = \frac{1 }{ \alpha - \dt } \left( R\bigl( u_h^0, u_h^1 \bigr) + \beta a \bigl( u_h^{1}, u_h^{1} \bigr) + \beta a \bigl( u_h^{0}, u_h^{0} \bigr) + \sum_{k=1}^p \bigl| g^k \bigr|^2 \dt \right).
\end{eqnarray*}

% {\sl On peut aussi obtenir cette autre inegalite:
% $$\sum_{p=0}^n \left| {u_h^{p+1} - u_h^{p} \over \dt} \right|^2 \le \sum_{p=1}^n {k_p \over \left( 1 
%- {\dt \over \alpha} \right)^{n-p+1} } + \left| {u_h^{1} - u_h^{0} \over \dt} \right|^2 {1 \over \left( 1
% - {\dt \over \alpha} \right)^n }$$
% with
% $$k_p = {1 \over \alpha} \left( R\bigl( u_h^0, u_h^1 \bigr) + + {1 \over 4} a \bigl( u_h^{1}, u_h^{1} 
% \bigr) + {1 \over 4} a \bigl( u_h^{0}, u_h^{0} \bigr) + \sum_{k=1}^p \bigl| g^k \bigr|^2 \dt \right)$$
% mais la majoration du membre de droite n'est pas aussi claire (bien que vraie quand meme).}

Since $f \in L^2 (0,T; H )$, we infer that the right hand side of (\ref{eq2}) remains bounded by a constant which depends only on the data $(f, u_0, v_0)$.
\end{proof}

We may observe that the lack of stability is due to the terms $(1-2 \beta) a(u_h^n, u_h^{n+1})$ and $( 1-2 \beta) a(u_h^{n-1}, u_h^n)$. For the case $\beta = 1/2$,  this difficulty does not occur  and  we obtain an unconditional stability result:

\bigskip

\begin{proposition} \label{stab2}  Let $\beta = 1/2$. 
Then there exists a constant depending only on the data, $C( f, u_0, v_0)$, such that for all $h>0$ and for all $N \ge 1$
\begin{eqnarray*}
 \left| \frac{u_h^{n+1} - u_h^{n} }{ \dt } \right|^2 + \frac{1 }{ 2} a \bigl( u_h^{n}, u_h^{n} \bigr) 
 + \frac{1 }{ 2} a \bigl( u_h^{n+1}, u_h^{n+1} \bigr)\le C (f, u_0, v_0)
 \end{eqnarray*}
for all $n \in \{ 1, \dots, N-1 \}$, where  $(u_h^{n+1})_{1 \le n \le N-1}$ are the solutions of problems $(P_{h \beta}^{n+1})_{1 \le n \le N-1}$. 
\end{proposition}

We define now an approximate solution $u_{h,N}^\beta$ ($\beta \in [0, 1/2]$) of problem $(P)$ by a linear interpolation of the solutions $u_h^{n+1}$ of $(P_{h \beta}^{n+1})$. More precisely, for all $h>0$ and $N \ge 1$ 
\begin{eqnarray*}
u_{h,N}^\beta (x,t)= u_h^{n} \frac{(n+1) \dt -t }{ \dt} + u_h^{n+1} \frac{t - n \dt }{ \dt}, 
\end{eqnarray*}
for all $t \in \bigl[ n \dt ,  (n+1) \dt \bigr]$, $0 \le n \le N-1$.  

Let $\alpha \in (0,1)$ and $N_h$ be defined by condition (\ref{eq1}) if $\beta \in [0, 1/2)$, otherwise let $N_h =1$ for all $h>0$.  The previous stability results imply that there exists a subsequence, still denoted $(u_{h,N}^\beta)_{h >0, N \ge N_h}$, and $u \in W=\bigl\{ w \in L^{\infty} (0,T; V), w_t \in L^{\infty} (0,T; H)  \bigr\}$ such that
\begin{eqnarray*}
\begin{array}{l}
 u_{h,N}^\beta \rightharpoonup u \quad \hbox{\rm weakly* in $L^{\infty} (0,T; V)$,} \\
\displaystyle  \frac{\partial u_{h,N}^\beta }{ \partial t} \rightharpoonup \frac{\partial u }{ \partial t} \quad \hbox{\rm weakly* in $L^{\infty} (0,T; H)$.} 
 \end{array} 
 \end{eqnarray*}
With Simon's lemma (\cite{Si}) we infer that $W$ is compactly embedded in $C^0 \bigl( [0, T]; H^1 (0,L) \bigr)$ and  we know also that  $W \subset C^{0,1/2} \bigl( [0,L] \times [0,T] \bigr)$ (\cite{MS-Ber}). It follows that, possibly extracting another subsequence, we have 
\begin{eqnarray*}
u_{h,N}^\beta \to u \quad \hbox{\rm strongly  in $C^0 \bigl( [0, T]; H^1 (0,L) \bigr)$,} 
\end{eqnarray*}
and thus $u$ belongs to $L^2(0,T; K)$ and $u(\cdot , 0)=u_0$.

\bigskip

Let us prove now that $u$ is a  solution of problem $(P)$.

\bigskip

\begin{theorem} Let $\beta \in [0, 1/2]$ and let $N_h$ be defined by condition (\ref{eq1})  if $\beta \not=1/2$, otherwise $N_h =1$ for all $h>0$. The sequence of approximate solutions $(u_{h,N}^\beta)_{h > 0, N \ge N_h}$ admits a subsequence which converges weakly*  in $W$ to a solution of problem $(P)$.
\end{theorem}

\begin{proof} We consider now the converging subsequence of $(u_{h,N}^\beta )_{h>0, N \ge N_h} $, still denoted $(u_{h,N}^\beta)_{h>0, N \ge N_h} $. 
Let $\tilde w \in {\mathcal H} \cap L^2 (0,T; K)$ such that  $\tilde w(\cdot, T) = u(\cdot , T )$. We will prove that 
\begin{eqnarray*}
\begin{array}{l}
\displaystyle - \int_0^T \bigl( u_t (\cdot, t), \tilde w_t (\cdot, t) -u_t (\cdot, t) \bigr) \, dt +  \int_0^T a \bigl( u (\cdot, t), \tilde w (\cdot, t) -u (\cdot, t) \bigr) \, dt  \\
\displaystyle  \ge \bigl(v_0,\tilde  w(\cdot, 0) - u_0 \bigr) + \int_0^T \bigl( f(\cdot, t) , \tilde w (\cdot, t) -u (\cdot, t) \bigr) .
 \end{array}
 \end{eqnarray*}
As a first step we have to construct a well-suited test-function $w_h^n$.

\bigskip

 Let $\varepsilon \in (0,T/2)$ and $\phi$ be a $C^{\infty}$-function such that
\begin{eqnarray*}
\left\{
\begin{array}{l}
 0 \le \phi (t) \le 1 \quad \forall t \in [0,T], \\
 \phi (t) =0 \quad \forall t \in [T- 3 \varepsilon/2, T], \quad \phi(t) =1 \quad \forall t \in [0, T- 2 \varepsilon]. 
 \end{array} \right.
 \end{eqnarray*}
 We denote $w = (1- \phi) u + \phi \tilde w$. Since $K$ is convex, we have immediately $ w \in {\mathcal H} \cap L^2 (0,T; K)$ and $w(\cdot, t) = u( \cdot, t)$ for all $t \in [T- 3 \varepsilon /2, T]$.

Let $\displaystyle \eta \in (0, \varepsilon / 2 )$ and $\mu \in (0,1)$. Following the same ideas as in \cite{MS-Ber}  we define $w^{\eta, \mu}$ by
\begin{equation} \label{eq3}
w^{\eta, \mu} (\cdot, t)= 
\displaystyle u(\cdot ,t) + \frac{1 }{ \eta} \int_t^{t+ \eta } \bigl( (1- \mu) w(\cdot ,s) - u( \cdot ,s)  \bigr) \, ds \quad \forall t \in [0,  T- \varepsilon /2 ].  
\end{equation}
Since $u \in W$ and  $w \in {\mathcal H}$, we have immediately $w^{\eta, \mu} - u  \in C^0 \bigl( [0,T]; V \bigr)$,  $w_t^{\eta, \mu} \in L^{2} \bigl( 0,T; H \bigr)$ and  $w^{\eta, \mu} \in L^{\infty} \bigl( 0,T; V \bigr) \cap C^0 \bigl( [0,T]; H^1 (0,L) \bigr)$. Moreover we can choose $\eta$ such that $w^{\eta, \mu}$ satisfies strictly the constraint. More precisely, for all $t \in [0, T- \varepsilon /2 ]$ and for all $x \in [0,L]$  we have
\begin{eqnarray*}
w^{\eta, \mu} (x,t) =\frac{1 }{ \eta} \int_t^{t + \eta} (1- \mu) w(x,s) \, ds + u(x,t) - \frac{1 }{ \eta} \int_t^{t+ \eta} u(x,s) \, ds.
\end{eqnarray*}
The first term of the right hand side belongs to $\bigl[(1-\mu) g_1(x), (1- \mu) g_2 (x) \bigr]$ with the convention that $(1-\mu) g_i(x) =g_i(x) $ ($i=1,2$) if $g_i(x) \in \{ + \infty, - \infty\}$, and recalling that $u \in C^{0, 1/2} \bigl( [0,L] \times [0,T] \bigr)$ we have
\begin{eqnarray*}
 \left| u(x,t) - \frac{1 }{ \eta} \int_t^{t+ \eta} u(x,s) \, ds \right| \le \frac{1 }{ \eta} \int_t^{t+ \eta} \bigl| u(x,t) - u(x,s) \bigr| \, ds \le \frac{2 C_0 \sqrt{\eta} }{ 3} 
\end{eqnarray*}
where $C_0$ is the H\"older continuity coefficient of $u$.

Thus, choosing $\eta$ such that
\begin{equation} \label{eq4}
 \frac{2 C_0 \sqrt{\eta} }{ 3} \le \frac{\mu }{ 2} g 
\end{equation}
ensures that
\begin{equation} \label{eq5}
 g_1(x)  + \frac{\mu }{ 2}  g    \le w^{\eta, \mu}(x,t) \le g_2(x) - \frac{\mu }{ 2}  g 
\end{equation}
 for all $ t \in [0, T- \varepsilon/2]$ and for all $ x \in [0,L]$, with the convention that 
 $\displaystyle g_i (x) \pm \frac{\mu }{ 2} g = g_i (x)$ ($i=1,2$)  if 
 $g_i (x) \in \{ + \infty, - \infty \}$.
\bigskip

Now, we assume that $\dt <\displaystyle \frac{ \varepsilon }{ 2 }$ and, for all $n \in \{ 1, \dots, N-1 \}$,  we define the test-function $w_h^n$ by
\begin{equation} \label{eq6}
w_h^n= 
\left\{
\begin{array}{l}
 u_h^{n+1} + Q_h \bigl( w^{\eta, \mu} (\cdot, n \dt) - u(\cdot, n \dt) \bigr) \quad \hbox{\rm if $n \dt \le T- \varepsilon$,} \\
u_h^{n+1} \quad \hbox{\rm if $n \dt > T- \varepsilon$.} 
\end{array}
\right. 
\end{equation}
 We have to check that $w_h^n$ belongs to $K_h$.

\begin{lemma} There exist $h_1 >0$ and $N_h' \ge N_h$ such that, for all $h \in (0, h_1)$ and for all $N \ge N_h'$, we have
\begin{eqnarray*}
w_h^n \in K_h \quad \forall n  \in \{1, \dots, N-1\}.
\end{eqnarray*}
\end{lemma}

\begin{proof} First of all it is clear that $w_h^n \in V_h$ and $w_h^n \in K_h$ if $n \dt > T - \varepsilon$.
Otherwise, when $n \dt \le T- \varepsilon$, we rewrite $w_h^n$ as follows:
\begin{equation} \label{eq7}
\begin{array}{l}
 w_h^n =  u_{h,N}^{\beta} \bigl(\cdot, (n+1) \dt \bigr) - u \bigl(\cdot, (n+1) \dt \bigr) + u \bigl(\cdot, (n+1) \dt \bigr) - u(\cdot,  n \dt) \\
 + w^{\eta, \mu} (\cdot, n \dt) + (Q_h -I) \bigl( w^{\eta, \mu} (\cdot, n \dt) - u(\cdot, n\dt) \bigr).
 \end{array}
 \end{equation}
We already know that $(u_{h,N}^{\beta})_{h>0, N \ge N_h}$ converges to $u$ strongly in $C^0 \bigl( [0, T]; H^1 (0,L) \bigr) $ and $u \in C^{0, 1/2} \bigl( [0,L] \times [0,T] \bigr)$, thus
\begin{equation} \label{eq8}
\sup_{ x \in [0,L]} \bigl| u \bigl( x, (n+1) \dt \bigr) -  u \bigl( x, n \dt \bigr) \bigr| \le C_0 \sqrt{\dt} 
\end{equation}
and 
\begin{equation} \label{eq9}
\sup_{ x \in [0,L]} \bigl| u_{h,N}^{\beta} \bigl( x, (n+1) \dt \bigr) -  u \bigl( x, (n+1) \dt \bigr) \bigr| \le C_1 \| u_{h,N}^{\beta} -u \|_{C^0 \bigl( [0,T]; H^1 (0,L) \bigr)}  
\end{equation}
where $C_1$ is the norm of the canonical injection of $H^1(0,L)$ into $C^0\bigl([0,L] \bigr)$.
Moreover 
\begin{equation} \label{eq10}
\begin{array}{l}
 \sup_{x \in [0,L]}  \bigl| (Q_h -I) \bigl( w^{\eta, \mu} (x, n \dt) - u(x, n \dt) \bigr) \bigr|  \\
 \le C_1 \bigl\| (Q_h -I) \bigl( w^{\eta, \mu} ( \cdot , n \dt) - u ( \cdot, n \dt) \bigr) \bigr\|_{H^1(0,L) } 
 \le C_1 \gamma_h \bigl\| w^{\eta, \mu} -u \bigr\|_{L^{\infty} ( 0,T;V )}.
 \end{array} 
 \end{equation}
By choosing $h_1$ and $N_h' \ge N_h$ such that
\begin{eqnarray*}
\frac{ \mu }{ 2} g \ge C_0 \sqrt{\dt} + C_1 \| u_{h,N}^{\beta} -u \|_{C^0 \bigl( [0,T]; H^1(0,L) \bigr)} + C_1 \gamma_h \| w^{\eta, \mu}- u \|_{L^{\infty} ( 0,T;V )} 
\end{eqnarray*}
for all $h \in (0, h_1)$ and $\dt = \frac{T }{ N}$ with $N \ge N_h'$, relations (\ref{eq5}) and (\ref{eq7})-(\ref{eq10}) imply that $w_h^n (x)$ belongs to  $ \bigl[g_1 (x), g_2 (x) \bigr]$ for all $x \in [0,L]$ if $n \dt \le  T - \varepsilon $ which concludes the proof.
\end{proof}

Let us choose now $w_h = w_h^n$ in $(P_{h \beta}^{n+1})$, $1 \le n \le N-1$ and define $\displaystyle N'= \left\lfloor \frac{T- \varepsilon }{ \dt} \right\rfloor$. With a discrete integration we obtain
\begin{equation} \label{eq11}
\begin{array}{l}
\displaystyle  \left( \frac{u_h^1 - u_h^0 }{ \dt}, w_h^0 - u_h^1 \right) + \sum_{n=1}^{N'} \bigl( g^n, w_h^n - u_h^{n+1} \bigr) \dt \\
\displaystyle  \le \sum_{n=1}^{N'} a \bigl( \beta u_h^{n+1} + (1-2 \beta) u_h^n + \beta u_h^{n-1} , w_h^{n} - u_h^{n+1} \bigr) \dt 
\\
\displaystyle  - \sum_{n=1}^{N'+1} \left( \frac{u_h^{n} - u_h^{n-1} }{ \dt}, \frac{ \bigl( w_h^{n} - u_h^{n+1} \bigr) - \bigl( w_h^{n-1} - u_h^{n} \bigr) }{ \dt } \right) \dt 
 \end{array} 
 \end{equation}
and we have to pass to the limit in each term as $h$ and $\dt$ tend to zero.
Recalling (\ref{eq0}) we immediatly infer that
\begin{eqnarray*}
  \left( \frac{u_h^1 - u_h^0 }{ \dt}, w_h^0 - u_h^1 \right) \to \bigl( v_0, w^{\eta, \mu} (\cdot, 0) - u(\cdot, 0) \bigr).
\end{eqnarray*}
Then, we rewrite the second term as follows
\begin{eqnarray*}
\begin{array}{l}
\displaystyle  \sum_{n=1}^{N'} \bigl( g^n, w_h^n - u_h^{n+1} \bigr) \dt 
\displaystyle =  \sum_{n=1}^{N'} \bigl( f^n, w_h^n - u_h^{n+1} \bigr) \dt \\
 + \beta \sum_{n=1}^{N'+1} \bigl( f^n - f^{n-1}, (w^{n-1} - u^n ) - (w^n - u^{n+1})  \bigr) \dt 
 + \beta ( f^0 - f^{1}, w^{0} - u^1 ) \dt. 
 \end{array}
 \end{eqnarray*}
% variante de factorisation:
% $$\eqalign{ & \sum_{n=1}^{N'} \bigl( g^n, w_h^n - u_h^{n+1} \bigr) \dt 
% =  \sum_{n=1}^{N'} \bigl( f^n, w_h^n - u_h^{n+1} \bigr) \dt \cr
% & + \beta \sum_{n=2}^{N'} \bigl( f^n - f^{n-1}, (w^{n-1} - u^n ) - (w^n - u^{n+1} \bigr) \dt 
% + \beta ( f^{N'+1} - f^{N'}, w^{N'} - u^{N'+1} ) \dt + \beta ( f^0 - f^{1}, w^{1} - u^2 ) \dt. }$$
But
\begin{eqnarray*}
\begin{array}{l}
\displaystyle  \sum_{n=1}^{N'} \bigl( f^n, w_h^n - u_h^{n+1} \bigr) \dt  = \sum_{n=1}^{N'} \int_{n \dt}^{(n+1) \dt} \bigl( f( \cdot, s), (w^{\eta, \mu} -u ) (\cdot, s) \bigr) \, ds \\
\displaystyle  + \sum_{n=1}^{N'} \int_{n \dt}^{ (n+1) \dt} \bigl( f( \cdot, s), (w^{\eta, \mu} -u ) (\cdot, n \dt) -  (w^{\eta, \mu} -u ) (\cdot, s) \bigr) \, ds \\
\displaystyle  + \sum_{n=1}^{N'} \bigl( f^n, (Q_h - I) \bigl( (w^{\eta, \mu} -u ) (\cdot, n \dt) \bigr) \bigr) \dt.
 \end{array}
 \end{eqnarray*}
Observing that
\begin{eqnarray*}
\bigl| (Q_h -I) (w^{\eta, \mu} -u ) (\cdot, n \dt) \bigr| \le \gamma_h \| w^{\eta, \mu} -u \|_{L^{\infty} (0,T; V)}  
\end{eqnarray*}
for all $n \in \{ 1, \dots, N'\}$,  we obtain
\begin{eqnarray*}
\begin{array}{l}
\displaystyle  \left| \sum_{n=1}^{N'} \bigl( f^n, (Q_h - I) \bigl( (w^{\eta, \mu} -u ) (\cdot, n \dt) \bigr) \bigr) \dt \right| 
\displaystyle  \le \sum_{n=1}^{N'}  \gamma_h |f^n| \| w^{\eta, \mu} -u \|_{L^{\infty} (0,T; V) } \dt \\
\displaystyle  \le \sqrt{ T}  \gamma_h \| f \|_{L^2(0,T; H)} \| w^{\eta, \mu} -u \|_{L^{\infty} (0,T; V)} \to 0.
 \end{array}
 \end{eqnarray*}
Moreover, with the definition of $w^{\eta, \mu}$, we have
\begin{equation} \label{eq12}
\begin{array}{l}
\bigl\| (w^{\eta, \mu} -u ) (\cdot, n \dt) -  (w^{\eta, \mu} -u ) (\cdot, s) \bigr) \bigr\|_V \\
\displaystyle  \le \frac{1 }{ \eta} \int^s_{n \dt} \bigl\| \bigl( (1-\mu) w  -u \bigr) (\cdot, \sigma) \bigr\|_V \, d \sigma 
\displaystyle   + \int^{s + \eta}_{n \dt + \eta } \bigl\| \bigl( (1-\mu) w  -u \bigr) (\cdot, \sigma) \bigr\|_V \, d \sigma \\
\displaystyle  \le \frac{ 2 \sqrt{ s - n \dt} }{ \eta }  \| (1- \mu) w - u \|_{L^2(0,T; V)} 
 \end{array} 
 \end{equation}
for all $n \in \{0, \cdots, N' \}$ and $s \in [n \dt, (n+1) \dt]$.
 
\smallskip

If we denote by $C$ the norm of the canonical injection of $\bigl(V, \| \cdot \|_V \bigr)$ into $\bigl( H, | \cdot | \bigr)$, we get
\begin{eqnarray*}
\begin{array}{l}
\displaystyle  \left| \sum_{n=1}^{N'} \int_{n \dt}^{ (n+1) \dt} \bigl( f( \cdot, s), (w^{\eta, \mu} -u ) (\cdot, n \dt) -  (w^{\eta, \mu} -u ) (\cdot, s) \bigr) \, ds \right| \\
\displaystyle   \le \sqrt{\dt}  \frac{ 2 C  \sqrt{T}  }{ \eta} \| (1- \mu) w - u \|_{L^2(0,T; V)} \| f\|_{L^2(0,T; H)} \to 0.
 \end{array} 
 \end{eqnarray*}
Finally, since $f \in L^2(0,T;H)$, $w^{\eta, \mu} -u \in C^0 \bigl( [0,T]; V \bigr)\subset L^{\infty} (0,L; H)$ and $| T - \varepsilon - (N'+1) \dt| \le \dt$, we may conclude that 
\begin{eqnarray*}
\begin{array}{l}
\displaystyle \lim_{\dt \to 0} \sum_{n=1}^{N'} \int_{n \dt}^{(n+1) \dt} \bigl( f( \cdot, s), (w^{\eta, \mu} -u ) (\cdot, s) \bigr) \, ds \\
\displaystyle = \int_0^{T- \varepsilon} \bigl( f( \cdot, s), (w^{\eta, \mu} -u) (\cdot, s) \bigr) \, ds
\end{array}
\end{eqnarray*}
and
\begin{eqnarray*}
 \sum_{n=1}^{N'} \bigl( f^n, w_h^n - u_h^{n+1} \bigr) \dt  \to \int_0^{T- \varepsilon} \bigl( f( \cdot, s), (w^{\eta, \mu} -u) (\cdot, s) \bigr) \, ds .
\end{eqnarray*}

\bigskip

Moreover, relation (\ref{eq12}) implies that
\begin{equation} \label{eq13}
\bigl\| (w^{n-1} - u^{n}) - (w^{n} - u^{n+1}) \bigr\|_V \le \frac{ 2 \sqrt{\dt} }{ \eta } \bigl\| (1- \mu) w -u \bigr\|_{L^2(0,T; V)}  
\end{equation}
for all $n \in \{ 1, \dots, N' +1 \}$. It follows that 
\begin{eqnarray*}
\begin{array}{l}
\displaystyle \left| \sum_{n=1}^{N'+1} \bigl( f^n - f^{n-1}, (w^{n-1} - u^n ) - (w^n - u^{n+1})  \bigr) \dt  \right| \\
\displaystyle  \le \sqrt{\dt} \frac{ 4 C \sqrt{T} }{ \eta } \bigl\| (1- \mu) w -u \bigr\|_{L^2(0,T; V)} \|f\|_{L^2(0,T;H)} . 
 \end{array}
 \end{eqnarray*}
Finally we observe that 
\begin{eqnarray*}
 \| w_h^n - u_h^{n+1} \|_V \le \| w^{\eta, \mu} -u \|_{L^{\infty} (0,T; V)} \quad \forall n \in \{ 0, \dots, N'\}.
\end{eqnarray*}
Hence 
\begin{eqnarray*}
\bigl|   ( f^0 - f^{1}, w^{0} - u^1 ) \dt \bigr| \le 2  C  \sqrt{\dt} \| w^{\eta, \mu} -u \|_{L^{\infty} (0,T; V)} \| f\|_{L^2(0,T;H)}
\end{eqnarray*}
and we may conclude that 
\begin{eqnarray*}
 \sum_{n=1}^{N'} \bigl( g^n, w_h^n - u_h^{n+1} \bigr) \dt  \to \int_0^{T- \varepsilon} \bigl( f( \cdot, s), (w^{\eta, \mu} -u) (\cdot, s) \bigr) \, ds .
\end{eqnarray*}

\bigskip

Let us study now the convergence of the first term of the right hand side of (\ref{eq11}). First we rewrite it as follows:
\begin{equation} \label{eq14}
\begin{array}{l}
\displaystyle  \sum_{n=1}^{N'} a \bigl( \beta u_h^{n+1} + (1-2 \beta) u_h^n  + \beta u_h^{n-1} , w_h^{n} - u_h^{n+1} \bigr) \dt \\
\displaystyle  = a \left( \frac{ 1 -2 \beta }{ 2} u_h^1 + \beta u_h^0  , w_h^0 -u_h^1 \right) \dt \\
 \displaystyle  + \sum_{n=1}^{N'+1} a \left( \beta u_h^{n-1} + \frac{ 1-2 \beta }{ 2} u_h^n , ( w_h^n -u_h^{n+1} ) - ( w_h^{n -1} -u_h^{n} ) \right) \dt \\
\displaystyle  + \sum_{n=1}^{N'} a \left( \frac{u_h^{n+1} + u_h^{n} }{ 2},  w_h^n -u_h^{n+1} \right) \dt .
 \end{array}
 \end{equation} 
% Variante de factorisation: 
% $$\eqalign{ & \sum_{n=1}^{N'} a \bigl( \beta u_h^{n+1} + (1-2 \beta) u_h^n  + \beta u_h^{n-1} , w_h^{n} % - u_h^{n+1} \bigr) \dt = a \left( { 1 -2 \beta \over 2} u_h^1 + \beta u_h^0  , w_h^0 -u_h^1 \right) \dt % \cr
% & - a \left( { 1-2 \beta  \over 2} u_h^{N'+1} + \beta  u_h^{N'} , w_h^{N'} -u_h^{N'+1} \right) \dt 
% + \sum_{n=1}^{N'} a \left( \beta u_h^{n-1} + { 1-2 \beta \over 2} u_h^n , ( w_h^n -u_h^{n+1} ) 
% - ( w_h^{n -1} -u_h^{n} ) \right) \dt \cr
% & + \sum_{n=1}^{N'} a \left( {u_h^{n+1} + u_h^{n} \over 2},  w_h^n -u_h^{n+1} \right) \dt .} % \leqno(14)$$
 With the propositions \ref{stab1} and \ref{stab2} we know that $\bigl( \| u_{h, N}^{\beta} \|_{L^{\infty} (0,T; V)} \bigr)_{h>0, N \ge N_h}$ is bounded independently of $h$ and $\dt$, thus 
\begin{eqnarray*}
\displaystyle  \left| a \left( \frac{ 1 -2 \beta }{ 2} u_h^1 + \beta u_h^0 , w_h^0 -u_h^1 \right) \right| 
\displaystyle  \le \frac{1 }{ 2}  \| u_{h, N}^{\beta} \|_{L^{\infty} (0,T; V)}  \| w^{\eta, \mu} -u \|_{L^{\infty} (0,T; V)} 
 \end{eqnarray*}
and, with (\ref{eq13})
\begin{eqnarray*}
\begin{array}{l} 
\displaystyle  \left|  a \left( \beta u_h^{n-1} + \frac{ 1-2 \beta }{ 2} u_h^n , ( w_h^n -u_h^{n+1} ) - ( w_h^{n -1} -u_h^{n} ) \right)  \right| \\
\displaystyle  \le \frac{ \sqrt{\dt} }{ \eta}  \| u_{h,N}^{\beta} \|_{L^{\infty} (0,T; V)} \bigl\| (1- \mu) w - u \bigr\|_{L^{2} (0,T; V)}   .
 \end{array} 
 \end{eqnarray*}

\smallskip

Finally, we rewrite the last term of (\ref{eq14}) as follows:
\begin{eqnarray*}
\begin{array}{l}
\displaystyle  \sum_{n=1}^{N'} a \left( \frac{u_h^{n+1} + u_h^{n} }{ 2},  w_h^n -u_h^{n+1} \right) \dt 
\displaystyle = \sum_{n=1}^{N'} \int_{n \dt}^{(n+1) \dt} a \bigl( u_{h,N}^{\beta} (\cdot, s), w_h^n - u_h^{n+1} \bigr) \, ds \\
\displaystyle  =  - \int_{N' \dt}^{T- \varepsilon} a \bigl( u_{h,N}^{\beta} (\cdot, s) , Q_h ( w^{\eta, \mu} - u) (\cdot , s) \bigr) \, ds 
\displaystyle  - \int_0^{\dt} a \bigl( u_{h,N}^{\beta} (\cdot, s) , Q_h ( w^{\eta, \mu} - u) (\cdot , s) \bigr) \, ds \\ 
\displaystyle  + \sum_{n=1}^{N'} \int_{n \dt}^{(n+1) \dt} a \bigl( u_{h,N}^{\beta} (\cdot, s), Q_h \bigl( (w^{\eta, \mu} -u )(\cdot, n \dt ) \bigr) - Q_h \bigl( (w^{\eta, \mu} -u )(\cdot , s ) \bigr) \, ds \\
\displaystyle  + \int_0^{T- \varepsilon} a \bigl( u_{h,N}^{\beta} (\cdot, s) , Q_h ( w^{\eta, \mu} - u) ( \cdot, s ) \bigr) \, ds . 
 \end{array}
 \end{eqnarray*}

The two  first terms can be estimated by
\begin{eqnarray*}
  \| u_{h,N}^{\beta} \|_{L^{\infty} (0,T; V)} \| w^{\eta, \mu} - u \|_{L^{\infty} (0,T; V)} \dt 
\end{eqnarray*}
and, with estimate (\ref{eq12})  we have
\begin{eqnarray*}
\begin{array}{l}
\displaystyle  \left| \sum_{n=1}^{N'} \int_{n \dt}^{(n+1) \dt} a \bigl( u_{h,N}^{\beta} (\cdot, s), Q_h \bigl( (w^{\eta, \mu} -u )(\cdot, n \dt ) \bigr) - Q_h \bigl( (w^{\eta, \mu} -u )(\cdot , s) \bigr) \bigr) \, ds \right| \\
\displaystyle  \le \frac{2 T \sqrt{\dt} }{ \eta}  \| u_{h,N}^{\beta} \|_{L^{\infty} (0,T; V)} \bigl\| (1- \mu) w - u \bigr\|_{L^{2} (0,T; V)}. 
 \end{array} 
 \end{eqnarray*}

\smallskip

Finally, recalling that $Q_h(w^{\eta, \mu} -u)$ is the orthogonal projection of $w^{\eta, \mu} -u$ on $V_h$ respectively to the scalar product defined by $a$ on $V$, we obtain that
\begin{eqnarray*}
\displaystyle \int_0^{T- \varepsilon} a \bigl( u_{h,N}^{\beta} (\cdot, s) , Q_h ( w^{\eta, \mu} - u) ( \cdot, s ) \bigr) \, ds 
\displaystyle = \int_0^{T- \varepsilon} a \bigl( u_{h,N}^{\beta} (\cdot, s) , ( w^{\eta, \mu} - u) ( \cdot, s ) \bigr) \, ds
\end{eqnarray*}
and the weak convergence of  $u_h$  to $u$ in $L^2 (0,T; V)$ allows us to conclude that 
\begin{eqnarray*}
\displaystyle   \sum_{n=1}^{N'} a \bigl( \beta u_h^{n+1} + (1-2 \beta) u_h^n + \beta u_h^{n} ,  w_h^n -u_h^{n+1} \bigr) \dt 
\displaystyle   \to \int_0^{T - \varepsilon} a \bigl(u (\cdot, s ), ( w^{\eta, \mu} -u) (\cdot, s ) \bigr) \, ds.
\end{eqnarray*}

\bigskip

There remains now to study the convergence of the last term of (\ref{eq11}) i.e
\begin{equation} \label{eq15}
\sum_{n=1}^{N'+1} \left( \frac{u_h^{n} - u_h^{n-1} }{ \dt}, \frac{ ( w_h^{n} - u_h^{n+1} ) - ( w_h^{n-1} - u_h^{n}) }{ \dt } \right) \dt. 
\end{equation}

In order to simplify the notations, let us define
\begin{eqnarray*}
\psi_{\dt} (\cdot , t) = \frac{ (w^{\eta, \mu}-u) (\cdot, t + \dt) - (w^{\eta, \mu} -u) (\cdot, t ) }{ \dt } , \quad \forall t \in [0,  T- \varepsilon ].
\end{eqnarray*}

We rewrite (\ref{eq15}) as follows:
\begin{eqnarray*} 
\begin{array}{l}
\displaystyle  \sum_{n=1}^{N'+1} \left( \frac{u_h^{n} - u_h^{n-1} }{ \dt}, \frac{ ( w_h^{n} - u_h^{n+1} ) - ( w_h^{n-1} - u_h^{n}) }{ \dt } \right) \dt \\
\displaystyle  = - \left( \frac{u^{N'+1} - u^{N'} }{ \dt} , w^{N'}_h - u^{N'+1}_h \right) \\
\displaystyle  + 
\sum_{n=1}^{N'} \left( \frac{ u_h^{n} - u_h^{n-1} }{ \dt}, ( Q_h - I) \bigl( \psi_{\dt} \bigl(\cdot,  (n-1) \dt  \bigr)  \bigr) \right) \dt \\
\displaystyle  + \sum_{n=1}^{N'} \int_{(n-1) \dt}^{n \dt}  \left( \frac{ u_h^{n} - u_h^{n-1} }{ \dt},  \psi_{\dt} \bigl( \cdot,  (n-1) \dt  \bigr) - \psi_{\dt } (\cdot,  t) \right) \, dt \\
 \displaystyle  + \sum_{n=1}^{N'} \int_{(n-1) \dt}^{n \dt}  \left( \frac{ u_h^{n} - u_h^{n-1} }{ \dt},   \psi_{\dt } (\cdot,  t) \right) \, dt .
\end{array}
\end{eqnarray*}
The first term, which can be interpreted as a boundary term at $t=T$ for the discrete time integration, can be estimated by
\begin{eqnarray*}
\begin{array}{l}
\displaystyle  \left| \left( 
\frac{u^{N'+1} - u^{N'} }{ \dt} , w^{N'}_h - u^{N'+1}_h \right) \right|  
 \\
\displaystyle  \le \frac{ C }{ \eta} \max_{ 1 \le n \le N} \left| \frac{ u_h^{n} - u_h^{n-1} }{ \dt} \right| 
\int_{N' \dt}^{N' \dt + \eta} \bigl\| (1- \mu) w ( \cdot, s) -u ( \cdot, s)\bigr\|_V  \, ds. 
\end{array}
\end{eqnarray*}
But $N' \dt \ge T - 3 \varepsilon/2$,  thus $w( \cdot, s) = u(\cdot, s)$ for all $s \in [N' \dt, N' \dt + \eta]$ and 
\begin{eqnarray*}
\displaystyle \int_{N' \dt}^{N' \dt + \eta} \bigl\| \bigl( (1- \mu) w  -u \bigr) ( \cdot, s) \bigr\|_V  \, ds = \int_{N' \dt}^{N' \dt + \eta} \mu \bigl\| u (\cdot  , s) \bigr\|_V  \, ds 
\displaystyle  \le \mu \eta \| u \|_{L^{\infty} (0,T; V)}.
\end{eqnarray*}
Since propositions \ref{stab1} and \ref{stab2} imply that $\displaystyle  \max_{ 1 \le n \le N} \left| \frac{ u_h^{n} - u_h^{n-1} }{ \dt} \right|$ is bounded independently of $h$ and $ \dt$, we infer that there exists a constant $C'$ such that
\begin{eqnarray*}
 \left| \left( \frac{u^{N'+1}_h - u^{N'}_h  }{ \dt}, w^{N'}_h - u^{N'+1}_h \right) \right| \le 
  C' \mu \| u \|_{L^{\infty} (0,T; V)}  .
  \end{eqnarray*}

For the second term we perform the same kind of computation:
\begin{eqnarray*}
\begin{array}{l}
\displaystyle  \left| \sum_{n=1}^{N'} \left( \frac{ u_h^{n} - u_h^{n-1} }{ \dt}, ( Q_h - I) \bigl( \psi_{\dt} \bigl(\cdot,  (n-1) \dt  \bigr)  \bigr) \right) \dt \right| \\
\displaystyle  \le \sum_{n=1}^{N'} \gamma_h \left| \frac{ u_h^{n} - u_h^{n-1} }{ \dt} \right| \bigl\| (w^{\eta, \mu} -u ) 
(\cdot, n \dt) -  (w^{\eta, \mu} -u ) \bigl( \cdot, (n -1) \dt \bigr) \bigr\|_V \\
\displaystyle  \le \frac{\gamma_h }{ \eta } \max_{ 1 \le n \le N} \left| \frac{ u_h^{n} - u_h^{n-1} }{ \dt} \right| 
\sum_{n=1}^{N'} \left( \int_{(n-1) \dt}^{n \dt} \bigl\| \bigl( (1- \mu) w -u \bigr) (\cdot, s) \bigr\|_V \, ds  \right. \\
\displaystyle   + \left. \int_{(n-1) \dt+ \eta }^{n \dt + \eta } \bigl\| \bigl( (1- \mu) w -u \bigr) (\cdot, s) \bigr\|_V \, ds \right) \\
\displaystyle  \le 2 \frac{\gamma_h \sqrt{T} }{ \eta} \max_{ 1 \le n \le N} \left| \frac{ u_h^{n} - u_h^{n-1} }{ \dt} \right| \bigl\| (1- \mu) w -u \bigr\|_{L^2 ( 0,T; V)}.
 \end{array}
 \end{eqnarray*}
Recalling that  $(\gamma_h)_{h >0}$ converges to zero, we obtain 
\begin{eqnarray*}
 \sum_{n=1}^{N'} \left( \frac{ u_h^{n} - u_h^{n-1} }{ \dt}, ( Q_h - I) \bigl( \psi_{\dt} \bigl( \cdot, (n-1) \dt  \bigr)  \bigr) \right) \dt \to 0.
\end{eqnarray*}
In order to estimate the third term, we transform $\psi_{\dt} \bigl(\cdot,  (n-1) \dt  \bigr) - \psi_{\dt} (\cdot, t)$ as follows:
\begin{eqnarray*}
\begin{array}{l}
 \psi_{\dt} \bigl(\cdot,  (n-1) \dt  \bigr) - \psi_{\dt} (\cdot, t) \\
 \displaystyle = \frac{1 }{ \eta \dt} \left( \int_{(n-1) \dt}^{t} \Bigl( \bigl( (1- \mu) w -u \bigr) (\cdot, s + \dt) - \bigl( (1- \mu) w -u \bigr) (\cdot, s ) \Bigr) \, ds \right. \\
\displaystyle  \left.  + \int^{(n-1) \dt+ \eta }_{t+ \eta} \Bigl( \bigl( (1- \mu) w -u \bigr) (\cdot, s + \dt) - \bigl( (1- \mu) w -u \bigr) (\cdot, s ) \Bigr) \, ds \right) \\
\displaystyle   = \frac{1 }{ \eta \dt} \int_{(n-1) \dt}^{t} \int_{s}^{s + \dt} \bigl( (1- \mu) w_t -u_t \bigr) ( \cdot, \sigma) \, d \sigma \, ds \\
\displaystyle    + \frac{1 }{ \eta \dt} \int^{(n-1) \dt + \eta }_{t + \eta } \int_{s}^{s + \dt} \bigl( (1- \mu) w_t -u_t \bigr) ( \cdot, \sigma) \, d \sigma \, ds .
 \end{array}
 \end{eqnarray*}

Hence
\begin{eqnarray*}
 \bigl| \psi_{\dt} \bigl( \cdot, (n-1) \dt \bigr) - \psi_{\dt} (\cdot, t) \bigr| \le 
 \frac{ 2 \bigl( t - (n-1) \dt \bigr) }{ \eta \sqrt{\dt} } \bigl\| (1- \mu) w_t -u_t \bigr\|_{L^2 (0,T; H) } 
 \end{eqnarray*}
and 
\begin{eqnarray*}
\begin{array}{l}
\displaystyle   \left| \sum_{n=1}^{N'} \int_{(n-1) \dt}^{n \dt}  
 \displaystyle  \left( \frac{ u_h^{n} - u_h^{n-1} }{ \dt},  \psi_{\dt} \bigl(\cdot,  (n-1) \dt \bigr) - \psi_{\dt } (\cdot,  t) \right) \, dt \right| \\
\displaystyle   \le \sum_{n=1}^{N' } \frac{ \dt^2 }{ \eta \sqrt{\dt} }  
\displaystyle  \left| \frac{ u_h^{n} - u_h^{n-1} }{ \dt} \right| \bigl\| (1- \mu) w_t -u_t \bigr\|_{L^2 (0,T; H) } \\
\displaystyle  \le \frac{ T \sqrt{ \dt } }{ \eta} \max_{1 \le n \le N} \left| \frac{ u_h^{n} - u_h^{n-1} }{ \dt} \right| \bigl\| (1- \mu) w_t -u_t \bigr\|_{L^2 (0,T; H) } \to 0.
 \end{array}
 \end{eqnarray*}

\smallskip

Finally, we observe that
\begin{eqnarray*}
\begin{array}{l}
\displaystyle  \sum_{n=1}^{N'} \int_{(n-1)\dt}^{n\dt} \left( \frac{ u_h^{n} - u_h^{n-1} }{ \dt}, \psi_{\dt} ( \cdot, t) \right) \, dt
\displaystyle  = \int_0^{T-\varepsilon} \left( \frac{\partial u_{h,N}^{\beta} }{ \partial t} ( \cdot, t), \psi_{\dt} ( \cdot , t) \right) \, dt \\
\displaystyle  - \int_{N' \dt}^{T- \varepsilon} \left( \frac{u^{N'}_h - u^{N'-1}_h }{ \dt}, \psi_{\dt} (\cdot, t) \right) \, dt .
 \end{array}
 \end{eqnarray*}

Since $(1- \mu) w -u \in L^2(0,T; V)$, and 
\begin{eqnarray*}
\displaystyle  \psi_{\dt} (\cdot, t)  = \frac{1 }{ \eta \dt} \left( \int_{t+ \eta}^{t +  \eta + \dt} \bigl( (1- \mu) w  - u  \bigr) (\cdot, s) \, ds  - \int_{t }^{t  + \dt} \bigl( (1- \mu) w  - u  \bigr) (\cdot, s) \, ds \right) 
\end{eqnarray*}
 we obtain that
\begin{eqnarray*}
\psi_{\dt} (\cdot, t) \to_{\dt \to 0} \frac{ \bigl( (1- \mu) w -u \bigr) (\cdot, t + \eta) - \bigl( (1- \mu) w -u \bigr) (\cdot, t ) }{ \eta} 
\end{eqnarray*}
 strongly in $\displaystyle L^2\left(0, T- \varepsilon ; V \right)$.  Since $\displaystyle \frac{\partial u_{h,N}^{\beta} }{ \partial t}$ converges weakly to $\displaystyle \frac{\partial u }{ \partial t}$ in $L^2(0,T;V)$,  it follows that
\begin{eqnarray*}
\begin{array}{l}
\displaystyle  \int_0^{T-\varepsilon} \left( \frac{\partial u_{h,N}^{\beta} }{ \partial t} ( \cdot, t), \psi_{\dt} ( \cdot , t) \right) \, dt
 \to \\
\displaystyle  \int_0^{T- \varepsilon} \left( \frac{\partial u }{ \partial t} ( \cdot,t), \frac{ \bigl( (1- \mu) w -u \bigr) (\cdot, t + \eta ) - \bigl( (1- \mu) w -u \bigr) (  \cdot, t) }{ \eta }  \right) \, dt .
 \end{array}
 \end{eqnarray*}
 Moreover, for all $t \in [0, T- \varepsilon]$ 
\begin{eqnarray*}
\begin{array}{l}
\displaystyle  \bigl| \psi_{\dt} ( \cdot, t) \bigr| \le \frac{C }{ \dt} \bigl\| (w^{\eta, \mu} - u) (\cdot,  t + \dt) -  (w^{\eta, \mu} - u) (\cdot, t ) \bigr\|_V \\
\displaystyle  \le \frac{2 C }{ \eta \sqrt{\dt} }   \bigl\|  (1- \mu) w -u \bigr\|_{L^2(0,T;V)}.
 \end{array}
 \end{eqnarray*}
Thus
\begin{eqnarray*} 
\begin{array}{l}
\displaystyle \left| \int_{N' \dt}^{T- \varepsilon} \left( \frac{u^{N'}_h - u^{N'-1}_h }{ \dt}, \psi_{\dt} (\cdot, t) \right)  \, dt \right| \\
\displaystyle \le \frac{ 2 C \sqrt{\dt} }{ \eta}  \bigl\|  (1- \mu) w -u \bigr\|_{L^2(0,T;V)} \max_{1 \le n \le N} \left| \frac{ u_h^{n} - u_h^{n-1} }{ \dt} \right|.
\end{array}
\end{eqnarray*}

\bigskip

Taking into account the previous convergence results, we obtain 
\begin{eqnarray*}
\begin{array}{l}
\displaystyle  \int_0^{T- \varepsilon} \bigl(f (\cdot, t) , ( w^{\eta, \mu} -u) (\cdot, t)  \bigr) \, dt+ \bigl( v_0, (w^{\eta, \mu} -u) (\cdot , 0) \bigr) \\
\displaystyle  \le  \int_0^{T- \varepsilon} a \bigl( u (\cdot, t) , (w^{\eta, \mu} -u) (\cdot,t)  \bigr) \, dt \\
\displaystyle   - \int_0^{T- \varepsilon} 
\displaystyle  \left( \frac{\partial u }{ \partial t} ( \cdot, s) , \frac{ \bigl( (1- \mu) w -u \bigr) (\cdot, t + \eta) - \bigl( (1- \mu) w -u \bigr) (\cdot, t ) }{ \eta}  \right) \, dt \\
 + C' \mu \|u \|_{L^{\infty} (0,T;V) } 
\end{array}
\end{eqnarray*}
for all $\varepsilon \in (0,T/2) $, $\mu \in (0,1)$ and $\eta \in (0, \varepsilon/2)$ satisfying (\ref{eq4}).

\smallskip

Rewriting 
$ \bigl(  (1- \mu) w -u \bigr) ( \cdot, t + \eta) - \bigl( (1- \mu) w -u \bigr) ( \cdot, t) $
as
\begin{eqnarray*}
\displaystyle   \int_t^{t+ \eta} \bigl( (1- \mu) w_t -u_t \bigr) ( \cdot, \sigma ) \, d \sigma 
\end{eqnarray*}
 and recalling that $(1- \mu) w_t - u_t \in L^2 (0,T; H)$, we obtain  that
\begin{eqnarray*}
 \frac{1 }{ \eta} \int_t^{t+ \eta} \bigl( (1- \mu) w_t -u_t \bigr) (\sigma, \cdot) \, d \sigma \to_{\eta \to 0} (1- \mu) w_t - u_t 
\end{eqnarray*}
 strongly in $\displaystyle L^2 \left(0,T - \varepsilon  ; H \right)$. Similarly, since $(1- \mu) w -u \in L^2(0,T;V) \cap C^0 \bigl( [ 0,T]; H \bigr)$, we have
\begin{eqnarray*}
(w^{\eta, \mu} - u) (\cdot ,t) = \frac{1 }{ \eta} \int_t^{t + \eta} \bigl( (1- \mu) w - u \bigr) (\cdot ,s) \, ds \to_{\eta \to 0} \bigl( (1- \mu) w - u \bigr) 
\end{eqnarray*}
 strongly in $\displaystyle L^2 \left(0,T - \varepsilon  ; V \right)$, and 
\begin{eqnarray*}
 (w^{\eta, \mu} - u) (\cdot ,0) \to_{\eta \to 0} (1- \mu) w(\cdot, 0) - u(\cdot, 0) \quad \hbox{\rm strongly in $H$.}
\end{eqnarray*}
Thus, when $\eta $ tends to zero, we get
\begin{eqnarray*}
\begin{array}{l}
\displaystyle  \int_0^{T- \varepsilon} \bigl( f (\cdot,t), ( (1-\mu ) w -u )(\cdot,t)  \bigr) \, dt  + \bigl( v_0, ( (1-\mu ) w -u )(\cdot, 0)  \bigr) \\
\displaystyle \le  \int_0^{T- \varepsilon} a \bigl( u (\cdot, t) , ((1-\mu ) w -u) (\cdot,t)  \bigr) \, dt \\
\displaystyle - \int_0^{T- \varepsilon} \left( \frac{\partial u }{ \partial t} (\cdot, t) ,  ( (1- \mu) w_t -u_t) (\cdot,t)  \right) \, dt 
\displaystyle  + C' \mu \|u \|_{L^{\infty} (0,T;V) }.
\end{array}
\end{eqnarray*}
Finally, we can pass  to the limit when  $\mu$ and $\varepsilon$ tend to zero and, observing that $w-u= \phi (\tilde w -u)$, we may conclude the proof.
\end{proof}

\section{Finite element implementation in $(P_{h \beta}^{n+1})$} \label{sec:4}

We present now some  simulations when the contact with the stops takes place only at the right end of the
beam, i.e
\begin{eqnarray*}
\left\{
\begin{array}{l}
g_{1}\left( x\right) =-\infty, \quad g_{2}\left( x\right) =+\infty \quad \forall x\in [0,L), \\
g_{1}\left( L\right) =-g, \quad g_{2}\left( L\right) =g,
\end{array}
\right. \end{eqnarray*}
with $g$ a positive real number.
We use the well-known Hermite piecewise cubics as basis
functions for  the space discretization.
So  we  consider a partition of the interval $\left[ 0,L\right] $
into $J$ subintervals of length $h$, i.e. $x_{0}=0$, $x_{i}=ih$, ..., $%
x_{J}=L$. At each node $x_{i}$, we associate two Hermite piecewise cubics $%
\varphi _{2i-1}$ and $\varphi _{2i}$ ($i=1,...,J$) defined by
\begin{eqnarray*}
\begin{array}{l}
\varphi _{2i-1}\in P_{3}\mbox{, }\varphi _{2i-1}\left( x_{j}\right) =\delta
_{ij}\mbox{ and }\varphi _{2i-1}^{\prime }\left( x_{j}\right) =0 \quad
\mbox{
for }1\leq j\leq J\mbox{.} \\
\varphi _{2i}\in P_{3}\mbox{, }\varphi _{2i}^{\prime }\left( x_{j}\right)
=\delta _{ij}\mbox{ and }\varphi _{2i}\left( x_{j}\right) =0 \quad
\mbox{
for }1\leq j\leq J\mbox{.}
\end{array}
\end{eqnarray*}
Thus, if
\begin{eqnarray*}
u=\sum_{i=1}^{J}u_{2i-1}\varphi _{2i-1}+\sum_{i=1}^{J}u_{2i}\varphi _{2i},
\end{eqnarray*}
the coefficient $u_{2i-1}$ gives the  value of $u$ at
node $x_{i}$ and the coefficient $u_{2i}$ gives the value
of the derivative of $u$ at node $x_{i}$ , $1\leq i\leq J$. Hence, for any $u\in V$, the
interpolate $u_{h}$ is given by
\begin{eqnarray*}
u_{h}=\sum_{i=1}^{J}u\left( x_{i}\right) \varphi
_{2i-1}+\sum_{i=1}^{J}u^{\prime }\left( x_{i}\right) \varphi _{2i}.
\end{eqnarray*}
We  consider the following finite dimensional subspace
\begin{eqnarray*}
V_{h}={\rm span } \left\{ \varphi _{1},\varphi _{2}, \dots ,\varphi _{2J-1},\varphi
_{2J}\right\} \subset V.
\end{eqnarray*}

Let $N \ge 1$, $\beta \in [0, 1/2]$ and $(u_h^{n+1})_{0 \le n \le N-1}$ be the solutions of problems $(P_{h \beta}^{n+1})_{0 \le n \le N-1}$. Since $u_h^{n+1}$ belongs to $K_h= V_h \cap K$ we have
\begin{eqnarray*}
u_{h}^{n+1}=\sum_{i=1}^{2J}\overline{u}_{i}^{n+1}\varphi _{i} \  \hbox{\rm with } \ u_h^{n+1} (L)= \overline{u}_{2J-1}^{n+1} \in [-g, g]
\end{eqnarray*}
and  $\left( P_{h,\beta }^{n+1}\right) $ can be rewritten as follows:
\begin{eqnarray*}
\left\{
\begin{array}{l}
\hbox{\rm find $\displaystyle \overline{u}^{n+1} \in \overline K_{h}= \Er^{2J-2} \times [-g, g] \times \Er$ such that, for all $\overline{w}\in \overline K_{h}$} \\
\displaystyle \left( \mathcal{M}\left( \frac{\overline{u}^{n+1}-2\overline{u}^{n}+%
\overline{u}^{n-1}}{\Delta t^{2}}\right) ,\overline{w}-\overline{u}%
^{n+1}\right) \\
\displaystyle +\left( \mathcal{S}\left( \beta \overline{u}^{n+1}+\left(
1-2\beta \right) \overline{u}^{n}+\beta \overline{u}^{n-1}\right) ,\overline{%
w}-\overline{u}^{n+1}\right) \geq \left( G^{n},\overline{w}-\overline{u}%
^{n+1}\right) \end{array} \right.
\end{eqnarray*}
 where
\begin{eqnarray*}
G_{i}^{n}=(\beta f^{n+1}+\left( 1-2\beta \right) f^{n}+\beta f^{n-1},\varphi _{i}) \quad  i=1, \dots , 2J
\end{eqnarray*}
and  $\mathcal{M}$ and $\mathcal{S}$
are respectively the global mass and stiffness matrices.
The previous inequality is also equivalent to the differential inclusion:
\begin{eqnarray*}
\begin{array}{l}
\displaystyle \mathcal{M}\left( \frac{\overline{u}^{n+1}-2\overline{u}^{n}+\overline{u}%
^{n-1}}{\Delta t^{2}}\right) +\mathcal{S}\left( \beta \overline{u}%
^{n+1}+\left( 1-2\beta \right) \overline{u}^{n}+\beta \overline{u}%
^{n-1}\right)\\
\displaystyle +\partial \psi _{\overline K_{h}}\left( \overline{u}^{n+1}\right) \ni
G^{n} \end{array}
\end{eqnarray*}
which  can be rewritten as
\begin{equation} \label{inclusion2}
\left( \mathcal{M}+\Delta t^{2}\beta \mathcal{S}\right) \overline{u}%
^{n+1}+\Delta t^{2}\partial \psi _{\overline K_{h}}\left( \overline{u}^{n+1}\right)
\ni F^{n} \end{equation}
with $F^{n}=\left( 2\mathcal{M}-\Delta t^{2}(1-2\beta )\mathcal{S}\right)
\overline{u}^{n}-\left( \mathcal{M}+\Delta t^{2}\beta \mathcal{S}\right)
\overline{u}^{n-1}+\Delta t^{2}G^{n}$.

More precisely, if we consider a single beam element $\left[ x_{i},x_{i+1}\right] $ of
length $h$ ($i=0, \dots, J-1$), the elemental mass and stiffness matrices
are
\begin{eqnarray*}
\mathcal{M}_{e}=\frac{h}{420}\left(
\begin{array}{cccc}
156 & 22h & 54 & -13h \\
22h & 4h^{2} & 13h & -3h^{2} \\
54 & 13h & 156 & -22h \\
-13h & -3h^{2} & -22h & 4h^{2}
\end{array}
\right), \ % \end{eqnarray*}
% \begin{eqnarray*}
\mathcal{S}_{e}=\frac{2k^{2}}{h^{3}}\left(
\begin{array}{cccc}
6 & 3h & -6 & 3h \\
3h & 2h^{2} & -3h & h^{2} \\
-6 & -3h & 6 & -3h \\
3h & h^{2} & -3h & 2h^{2}
\end{array}
\right)
\end{eqnarray*}
and the global system (\ref{inclusion2}) is  obtained as an assembly of the
previous elemental matrices.

\bigskip

For the resolution of (\ref{inclusion2}) at each time step, we use the following
lemma with $\mathcal{A}=\left( \mathcal{M}+\Delta t^2 \beta \mathcal{S}%
\right) $, $\lambda=\Delta t^2$ and $f=F^n$

\begin{lemma} \cite{P93}, \cite{DY03}
Let $\mathcal{A}$ be a symmetric positive definite $2J\times 2J$ real matrix, $f \in \Er^{2J}$ and  $u^{\prime }$ be the solution of $\mathcal{A}u^{\prime }=f$.
Then, for all $\lambda>0$, the system \begin{eqnarray*}
\mathcal{A}u+\lambda \partial \psi _{\overline K_{h}}(u)\ni f,
\end{eqnarray*}
with $\overline K_{h}=\Er^{2J-2}\times \left[
-g,g\right] \times \Er$ admits an unique solution $u$ given by
\begin{eqnarray*}
\left\{
\begin{array}{ll}
u_{2J-1} & =P_{[-g,g]}(u_{2J-1}^{\prime }), \\
\left( \mathcal{A}u\right) _{i} & =f_{i},\qquad \mbox{for }i=1,...,2J-2,2J.
\end{array}
\right.
\end{eqnarray*}
\end{lemma}

We consider a steel pipe of length $L=1.501$ $m$, with an external diameter equal to $1$
$cm$ and a thickness equal to $0.5$ $mm$. Thus $\displaystyle k^{2}=\frac{EI}{\rho S}=282.84$ $ m^4. s^{-2}$
where $E=2\times 10^{11}$ $Pa$ is the Young's modulus, $\rho
=8 \times 10^3$ $kg/m^{3}$ is the material density, $S$ is the cross-section and $I$ the
cross-sectional moment of inertia of the pipe (see also \cite{P93} for a more detailed description  of the mechanical setting).
The vibration of the support is given by $\phi \left( t\right) =0.2\sin \left(
10t\right) $ for all $t \ge 0$, $g=0.1$ and the initial data are $u_0=0$, $v_0= -2 h$ i.e $\tilde u(\cdot , 0) = \tilde u_t (\cdot, 0) =0$ (at $t=0$ the beam is at rest).

\begin{center}
\begin{figure}[h]
\resizebox{0.75\textwidth}{!}{%
 \includegraphics{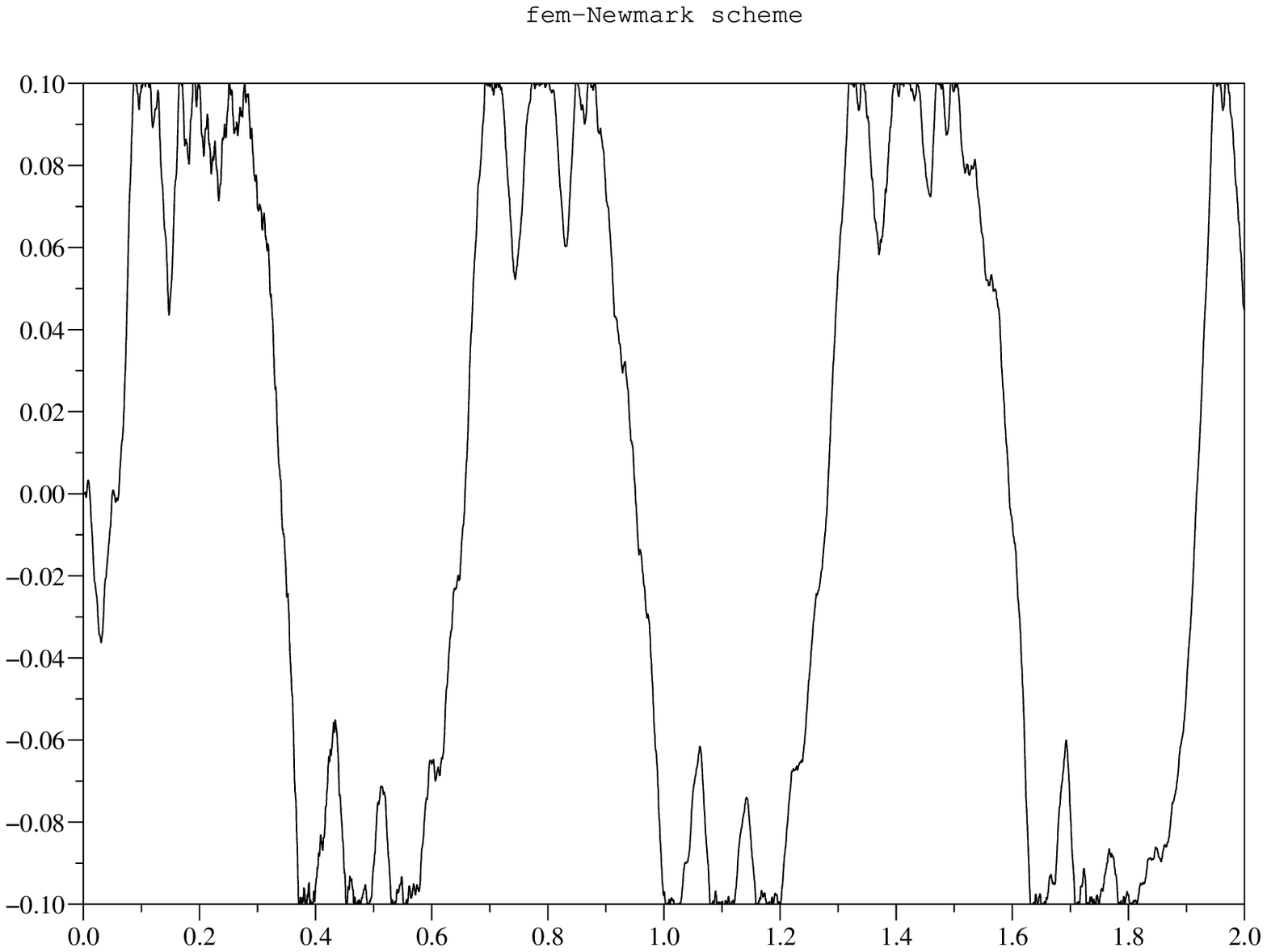}
}
%\vspace{5cm}       % Give the correct figure height in cm
\caption{$\Delta t=5 \times 10^{-5} s$ and $\beta =\frac12$}
\label{fig:2}       % Give a unique label
\end{figure}
\end{center}

\begin{center}
\begin{figure}[h]
\resizebox{0.75\textwidth}{!}{%
 \includegraphics{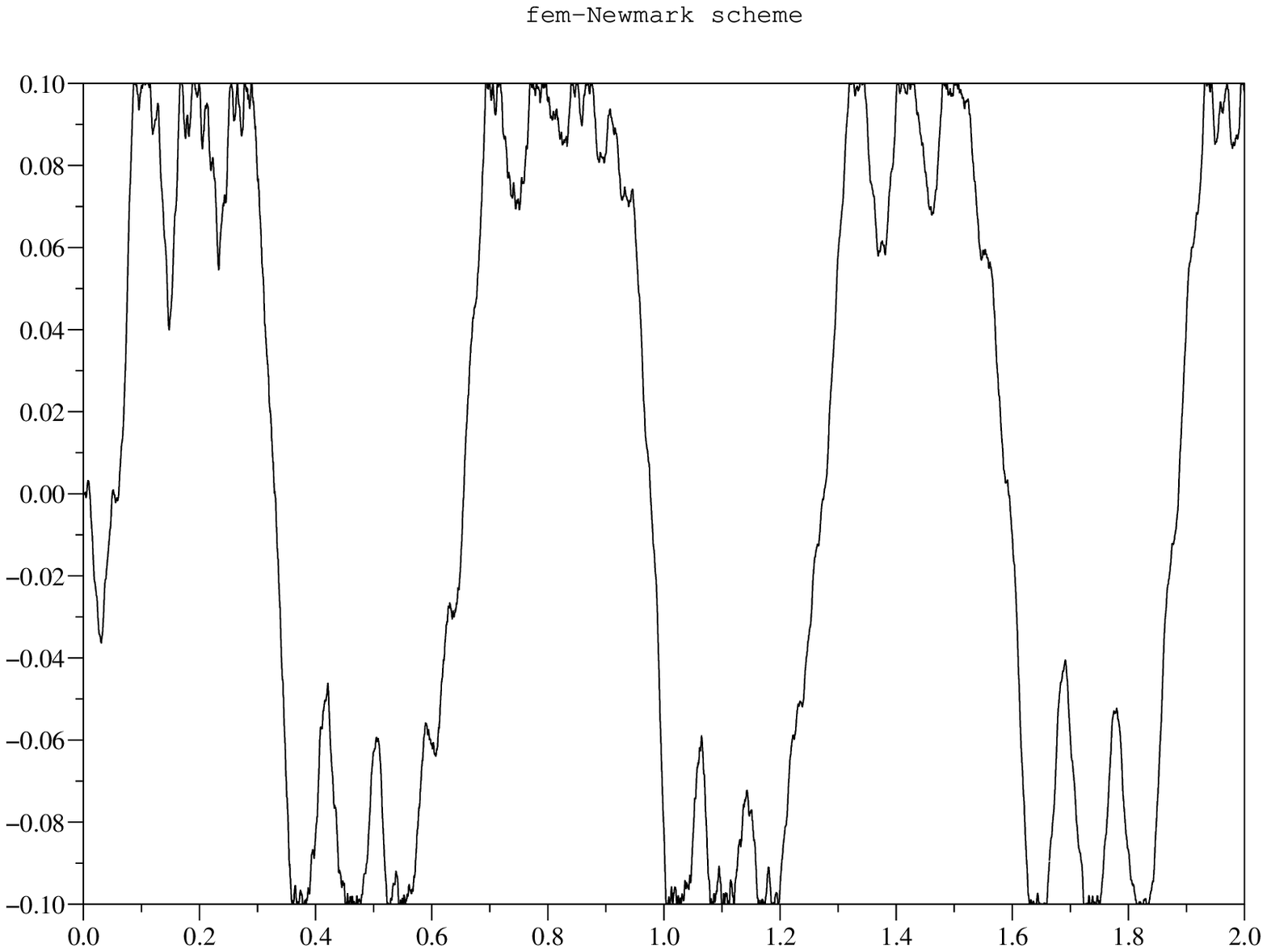}
}
%\vspace{5cm}       % Give the correct figure height in cm
\caption{$\Delta t=5\times 10^{-6} s$ and $\beta =\frac12$}
\label{fig:3}       % Give a unique label
\end{figure}
\end{center}

\begin{center}
\begin{figure}[h]
\resizebox{0.75\textwidth}{!}{%
 \includegraphics{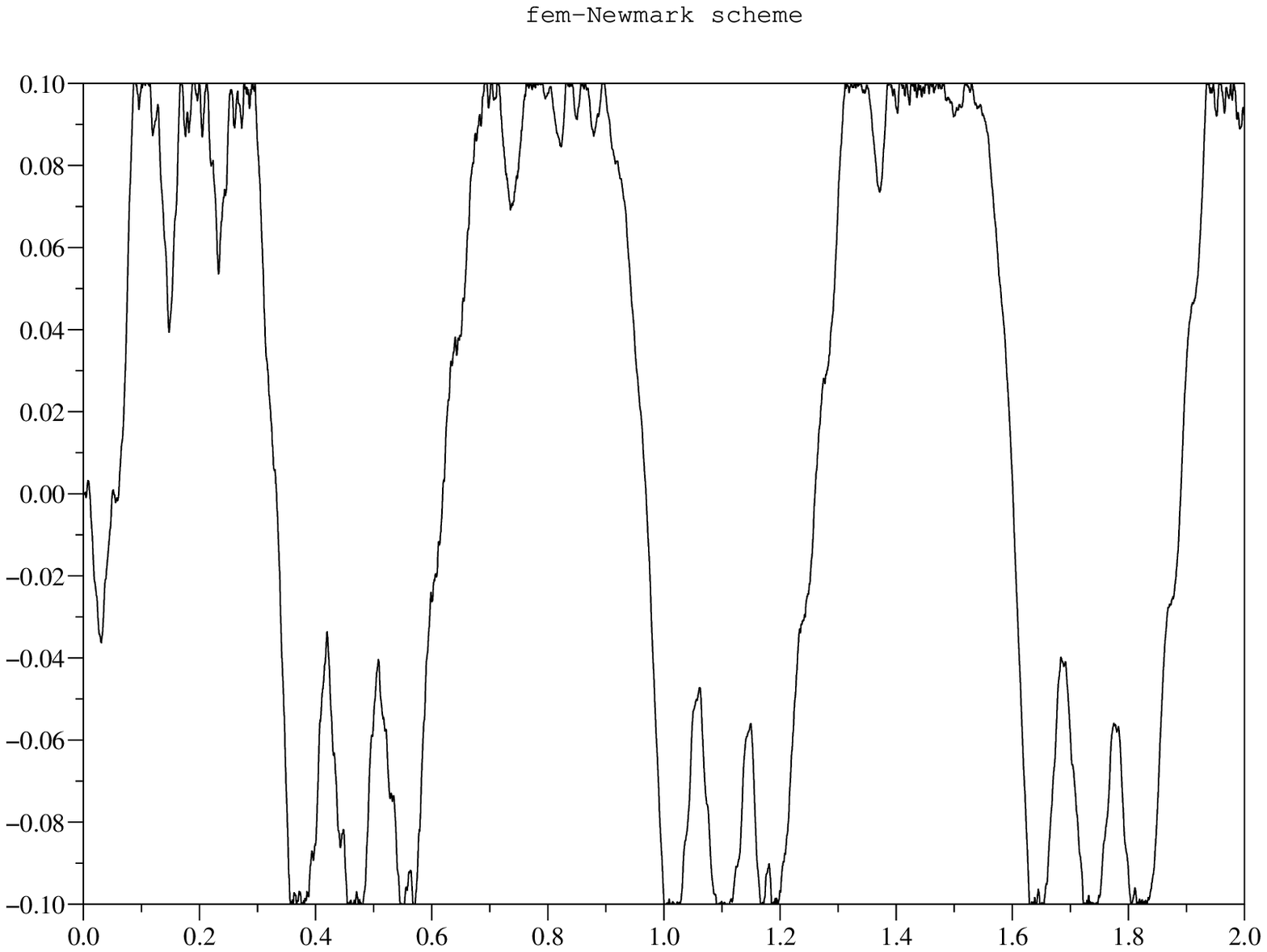}
}
%\vspace{5cm}       % Give the correct figure height in cm
\caption{$\Delta t=10^{-6} s$ and $\beta =\frac12$}
\label{fig:4}       % Give a unique label
\end{figure}
\end{center}

In the next figures we show the approximate motion of the impacting end of the beam.
The results given at figures \ref{fig:2}, \ref{fig:3} and \ref{fig:4} have been obtained with $\beta = 1/2$,
$J=19$ and $\Delta t = 5 \times 10^{-5} s$, $\Delta t = 10^{-5} s$ and
$\Delta t = 5 \times 10^{-6} s$ (let us recall that we have unconditional stability
for this value of $\beta$).

We can observe that the trajectories are almost the same at the beginning of
the time interval (up to the end of the first "contact period" i.e $0 \le t \le 0.2$) and remain quite similar afterwards even if the details of the impact
phenomenon are different. This is not surprising since vibrations with unilateral
constraints always lead to sensitivity to initial data.

The motion of the impacting end of the beam has also been computed by using the
normal compliance approximation of Signorini's conditions. In this case we have
to define the penalty parameter $\varepsilon$. Although the corresponding stiffness
$1/ \varepsilon$ has a physical meaning, the range of values usually chosen is quite
large:  $1/ \varepsilon = 10^{10}  N.m^{-1}$ in \cite{Ravn}, $1/\varepsilon = 5.5 \times 10^7
N.m^{-1}$ in \cite{Hurm} for instance. In the following results we consider $1/\varepsilon = 10^8 N.m^{-1}$
and we apply once again a Newmark's scheme with $J=19$ and $\Delta t = 5 \times 10^{-6} s$,
$\Delta t = 10^{-6} s$ and $\Delta t = 5 \times 10^{-7} s$ (see figures \ref{fig:5}, \ref{fig:6}, \ref{fig:7}). We should notice that we have to solve
now a  partial differential equation, thus we choose $\beta=1/4$ for which the
unconditional stability of Newmark's scheme holds.
\begin{center}
\begin{figure}[h]
\resizebox{0.75\textwidth}{!}{%
 \includegraphics{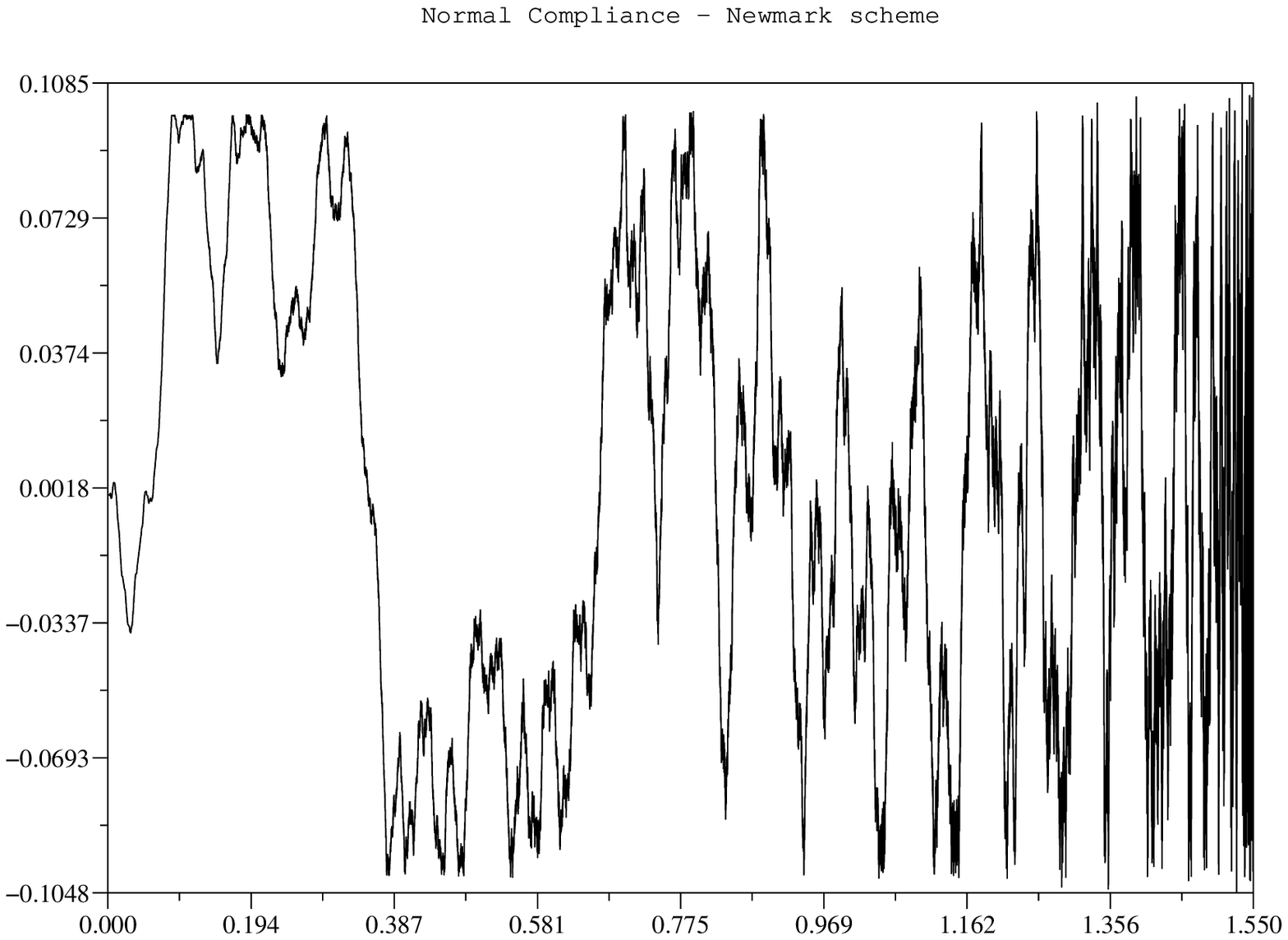}
}
%\vspace{5cm}       % Give the correct figure height in cm
\caption{$\Delta t=5 \times 10^{-6} s$ and $\beta =\frac14$}
\label{fig:5}       % Give a unique label
\end{figure}
\begin{figure}[h]
\resizebox{0.75\textwidth}{!}{%
 \includegraphics{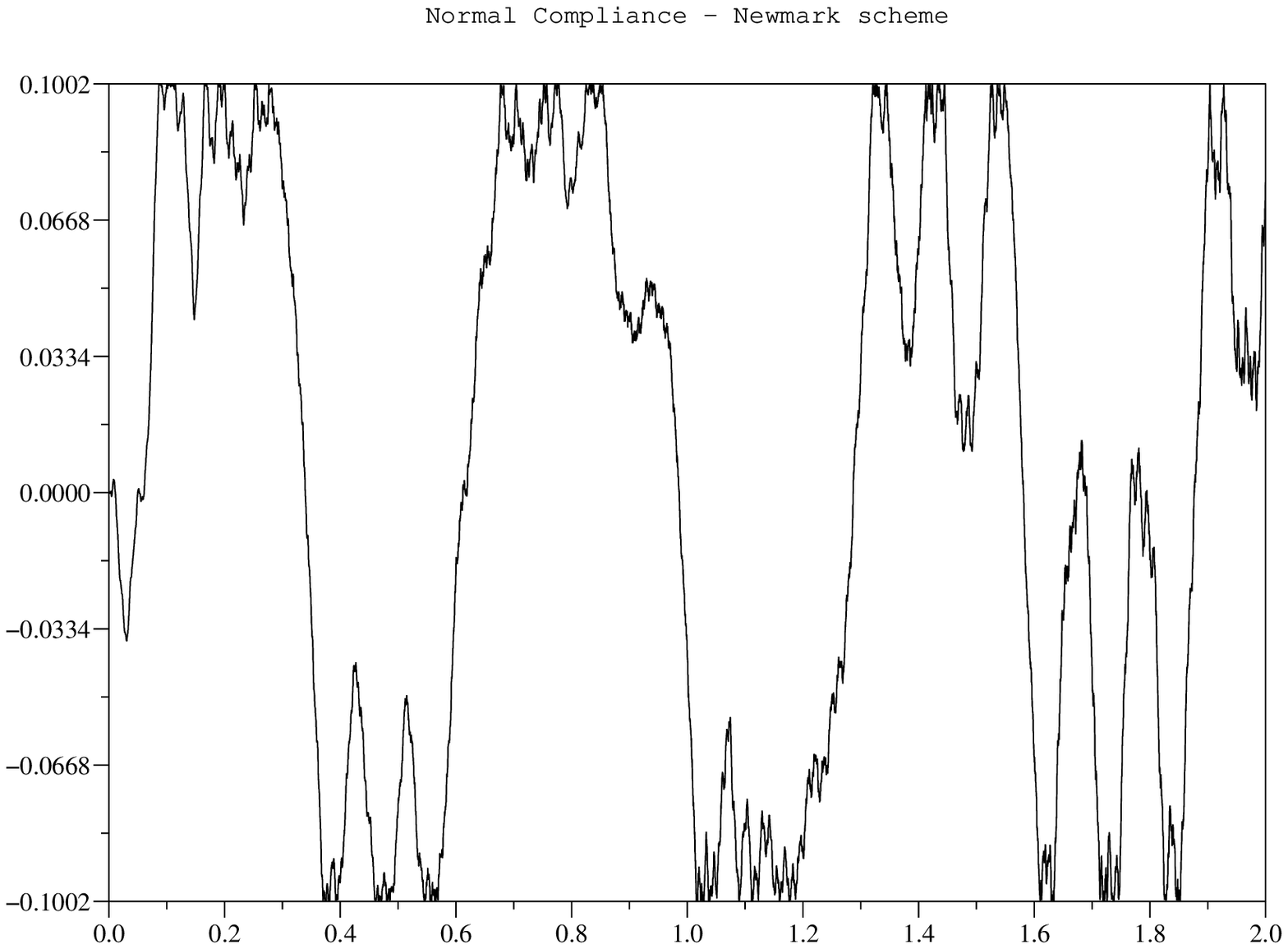}
}
%\vspace{5cm}       % Give the correct figure height in cm
\caption{$\Delta t=10^{-6} s$ and $\beta =\frac14$}
\label{fig:6}       % Give a unique label
\end{figure}
\begin{figure}[h]
\resizebox{0.75\textwidth}{!}{%
 \includegraphics{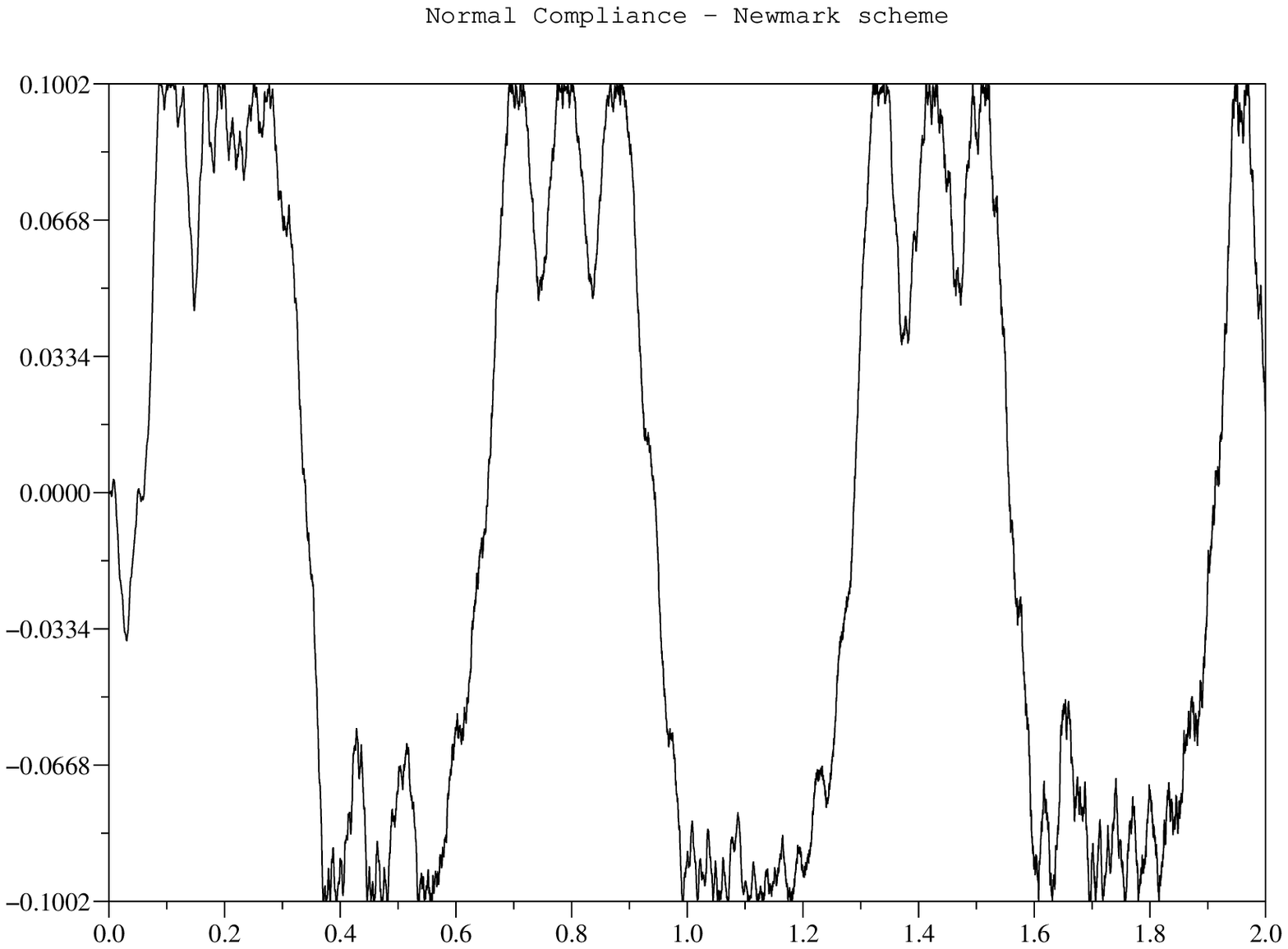}
}
%\vspace{5cm}       % Give the correct figure height in cm
\caption{$\Delta t=5 \times 10^{-7} s$ and $\beta =\frac14$}
\label{fig:7}       % Give a unique label
\end{figure}
\end{center}
Nevertheless we observe a kind
of numerical instability (figure \ref{fig:5}): spurious high frequencies appear during "contact periods"
and this phenomenon can be controlled only for very small time steps. Moreover, the non-penetration condition is violated by the trajectories computed with the
normal compliance approximation while it remains satisfied by the approximate
motions $u_{h,N }^{\beta}$. More precisely we have
\begin{eqnarray*}
\max_{0 \le t \le 2} \bigl| u_{app} (L,t) -g \bigr| = 2.25 \times 10^{-4} m \quad \hbox{\rm for
$\Delta t = 5 \times 10^{-7} s$}
\end{eqnarray*}
which is a rather coarse approximation of the unilateral constraint.

\begin{figure}[h]
\resizebox{0.75\textwidth}{!}{%
 \includegraphics{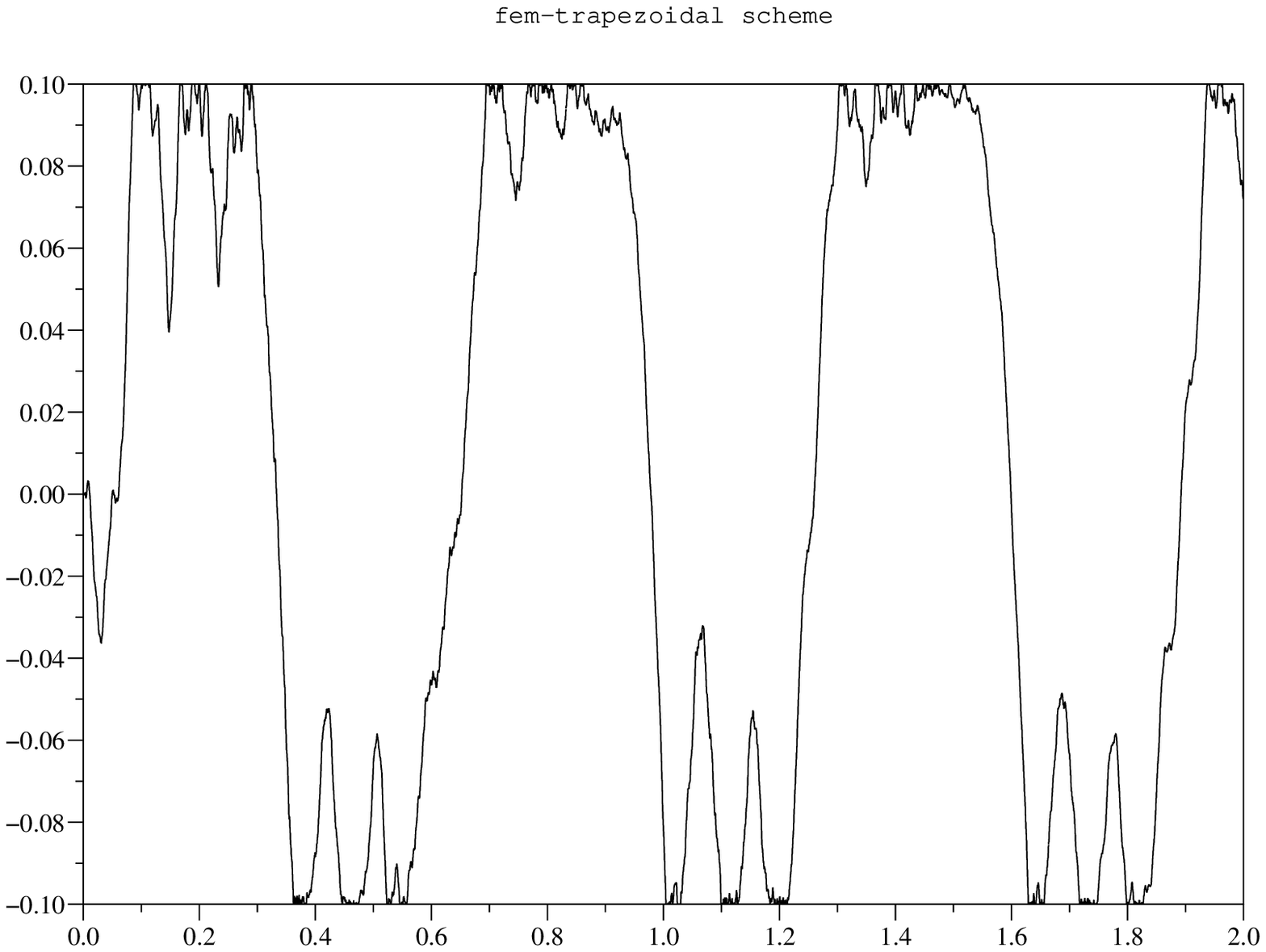}
}
%\vspace{5cm}       % Give the correct figure height in cm
\caption{$\Delta t=10^{-6} s$ and $\beta =\frac14$}
\label{fig:8}       % Give a unique label
\end{figure}

Finally let us point out that, for $J=19$ and $\beta = 1/4$, the stability condition (\ref{eq1})
leads to
\begin{eqnarray*}
\Delta t \le 2 \sqrt{\frac{2}{\tilde \kappa (h)} } = 3.3469 \times 10^{-6} s
\end{eqnarray*}
with $\tilde \kappa (h)$ given by (\ref{kappa(h)}) (see the Appendix). Thus, with $\beta = 1/4$
we can solve the discretized problems $(P_{h \beta}^{n+1})_{0 \le n \le N-1}$ with $J=19$ and
$\Delta t = 5 \times 10^{-6} s$ (figure \ref{fig:8}): the computation remains less expensive than with  the normal
compliance approximation (see the next table for a comparison of the CPU times).

\begin{eqnarray*}
\begin{tabular}{l|lllllll}
& Fig. 2 & Fig. 3 & Fig. 4 & Fig. 5 & Fig. 6 & Fig. 7 & Fig. 8 \\ \hline
CPU \\
Time & $18.24$ & $173.71$ & $870.45$ & $241.17$ & $1187.9$ & $2426.4$ & $859.38$%
\end{tabular}
\end{eqnarray*}

\begin{remark}
The computations has been performed on MAC G4 (1.1 Ghz) with Scilab, the scientific computing software developed by ENPC-INRIA. Other results are available at the following web address: http://www.univ-reunion.fr/~ydumont/beam-vibrations.htm
\end{remark}

\section{Appendix: estimate of $\protect\kappa \left( h\right) $.} \label{sec:5}

We consider the finite element space discretization described at section \ref{sec:4} i.e
$ V_h = {\rm span}\bigl\{ \varphi_1, \dots, \varphi_{2J} \bigr\}$
where $(\varphi_i)_{1 \le i \le 2J}$ is the Hermite piecewise cubics basis.
Thus, for all $u_h \in V_h$ we have
\begin{eqnarray*}
u_{h}=\sum_{i=1}^{J}u_h\left( x_{i}\right) \varphi
_{2i-1}+\sum_{i=1}^{J}u_h^{\prime }\left( x_{i}\right) \varphi _{2i}.
\end{eqnarray*}
In order to simplify the notations, we let
\begin{eqnarray*}
u_i = u_h\left(x_i\right) \quad u^{\prime}_i = u_h^{\prime}\left(x_i \right) \quad \hbox{\rm for all $i=0, \dots, J$.}
\end{eqnarray*}
We may observe that, since $u_h \in V_h \subset V$, we have $u_0=u'_0=0$.

Let us  compute first $(u_{h},u_{h})$. We have
\begin{eqnarray*}
\begin{array}{l}
\displaystyle \left( u_{h},u_{h}\right) = \frac{13}{35}\sum_{j=1}^{J}\left(
u_{j}^{2}+u_{j-1}^{2}\right) \Delta x+\frac{9}{35}\sum_{j=1}^{J}u_{j}u_{j-1}%
\Delta x \\
\displaystyle +\frac{13\Delta x}{210}\sum_{j=1}^{J}\left( u_{j}u_{j-1}^{\prime
}-u_{j-1}u_{j}^{\prime }\right) \Delta x 
+\frac{\left( \Delta x\right) ^{2}}{105}\sum_{j=1}^{J}\left( \left(
u_{j}^{\prime }\right) ^{2}+\left( u_{j-1}^{\prime }\right) ^{2}\right)
\Delta x  \\
\displaystyle -\frac{\left( \Delta x\right) ^{2}}{70}\sum_{j=1}^{J}u_{j}^{\prime
}u_{j-1}^{\prime }\Delta x-\frac{11\left( \Delta x\right) ^{2}}{105}%
u_{J}u_{J}^{\prime }.
\end{array}
\end{eqnarray*}

Then we use the following relations 
\begin{eqnarray*}
\begin{array}{l}
\displaystyle \frac{26}{420}\sum_{j=1}^{J}\Delta x\left( u_{j-1}^{\prime
}u_{j-1}-u_{j}^{\prime }u_{j}\right) \Delta x  =-\frac{26}{420}\left( \Delta
x\right) ^{2}u_{J}^{\prime }u_{J}, \\
\displaystyle \frac{26}{35}\sum_{j=1}^{J}u_{j}^{2}\Delta x-\frac{13}{35}u_{J}^{2}\Delta x 
=\frac{13}{35}\sum_{j=1}^{J}\left( u_{j}^{2}+u_{j-1}^{2}\right) \Delta x. 
\end{array}
\end{eqnarray*}
and we get
\begin{eqnarray*}
\begin{array}{l}
\displaystyle \left( u_{h},u_{h}\right)  =  \frac{1}{140}\sum_{j=1}^{J}\left( \frac{13}{3}%
\left( u_{j}+u_{j-1}\right) +\Delta x\left( u_{j-1}^{\prime }-u_{j}^{\prime
}\right) \right) ^{2}\Delta x 
 + \frac{1752}{7560}\sum_{j=1}^{J-1}\left(
u_{j}^{2}+u_{j-1}^{2}\right) \Delta x  \\ 
\displaystyle  +\frac{1}{180}\sum_{j=1}^{J}\left( u_{j}-u_{j-1}\right) ^{2}\Delta x 
\displaystyle  +\frac{%
\left( \Delta x\right) ^{2}}{420}\sum_{j=1}^{J-1}\left( \left( u_{j}^{\prime
}\right) ^{2}+\left( u_{j-1}^{\prime }\right) ^{2}\right) \Delta x+\frac{%
\left( \Delta x\right) ^{3}}{420}\left( u_{J-1}^{\prime }\right) ^{2}  
 \\
\displaystyle  +\frac{1}{420}\left( \Delta xu_{J}^{\prime }-9u_{J}\right) ^{2}\Delta x+%
\frac{1752}{7560}\Delta xu_{J-1}^{2}+\frac{294}{7560}u_{J}^{2}\Delta x. 
\end{array}
\end{eqnarray*}
Then, we observe that 
\begin{eqnarray*}
\begin{array}{l}
\displaystyle  \frac{1}{420}\left( \Delta xu_{J}^{\prime
}-9u_{J}\right) ^{2}\Delta x+\frac{294}{7560}u_{J}^{2}\Delta x \\
\displaystyle = \frac{1}{420}\left( \frac{18}{19}%
\Delta xu_{J}^{\prime }-\frac{19}{2}u_{J}\right) ^{2}\Delta x+\frac{51}{3024}u_{J}^{2}\Delta x+
\frac{37 \left(\Delta x\right) ^{2}}{420 \times 19^2} \left( u_{J}^{\prime }\right) ^{2}\Delta x
\end{array}
\end{eqnarray*}
and thus
\begin{eqnarray*}
\begin{array}{l}
\displaystyle  \left( u_{h},u_{h}\right)  = \frac{1}{140}\sum_{j=1}^{J}\left( \frac{13}{3}%
\left( u_{j}+u_{j-1}\right) +\Delta x\left( u_{j-1}^{\prime }-u_{j}^{\prime
}\right) \right) ^{2}\Delta x \\
\displaystyle  +\frac{1752}{7560}\sum_{j=1}^{J-1}\left(
u_{j}^{2}+u_{j-1}^{2}\right) \Delta x+\frac{1752}{7560}\Delta xu_{J-1}^{2} 
+\frac{1}{180}\sum_{j=1}^{J}\left( u_{j}-u_{j-1}\right) ^{2}\Delta x \\
\displaystyle  +\frac{%
\left( \Delta x\right) ^{2}}{420}\sum_{j=1}^{J-1}\left( \left( u_{j}^{\prime
}\right) ^{2}+\left( u_{j-1}^{\prime }\right) ^{2}\right) \Delta x 
\displaystyle  +\frac{%
\left( \Delta x\right) ^{2}}{420}\left( \left( u_{J-1}^{\prime }\right) ^{2}+%
\frac{37}{19^{2}}\left( u_{J}^{\prime }\right) ^{2}\right) \Delta x  \\
\displaystyle +\frac{1}{420}\left( \frac{18}{19}\Delta xu_{J}^{\prime }-\frac{19}{%
2}u_{J}\right) ^{2}\Delta x+\frac{51}{3024}u_{J}^{2}\Delta x.  
\end{array}
\end{eqnarray*}
We deduce the following inequality 
\begin{equation} \label{minoration u_h}
\begin{array}{l}
\displaystyle \left( u_{h},u_{h}\right)  \geq \frac{\left( \Delta x\right) ^{2}}{420}%
\frac{37}{19^{2}}\sum_{j=1}^{J}\left( \left( u_{j}^{\prime }\right)
^{2}+\left( u_{j-1}^{\prime }\right) ^{2}\right) \Delta x+\frac{1}{180}%
\sum_{j=1}^{J}\left( u_{j}-u_{j-1}\right) ^{2}\Delta x  \\
\displaystyle  \geq \frac{37}{420\times 19^{2}}\left( \left( \Delta x\right)
^{2}\sum_{j=1}^{J}\left( \left( u_{j}^{\prime }\right) ^{2}+\left(
u_{j-1}^{\prime }\right) ^{2}\right) \Delta x+\sum_{j=1}^{J}\left(
u_{j}-u_{j-1}\right) ^{2}\Delta x\right)  
\end{array} 
\end{equation}
Now we  compute $a\left( u_{h},u_{h}\right)$:
\begin{eqnarray*}
\begin{array}{l}
\displaystyle  a\left( u_{h},u_{h}\right)  = \frac{k^{2}}{\left( \Delta x\right) ^{4}}%
\left( 12\sum_{j=1}^{J}\left( u_{j}^{2}+u_{j-1}^{2}\right) \Delta
x-24\sum_{j=1}^{J}u_{j}u_{j-1}\Delta x\right)  \\
\displaystyle  +\frac{k^{2}}{\left( \Delta x\right) ^{4}}\left( 12\Delta
x\sum_{j=1}^{J}\left( u_{j-1}u_{j}^{\prime }-u_{j}u_{j-1}^{\prime }\right)
\Delta x+4\left( \Delta x\right) ^{2}\sum_{j=1}^{J}u_{j}^{\prime
}u_{j-1}^{\prime }\Delta x\right)  \\
\displaystyle  +\frac{k^{2}}{\left( \Delta x\right) ^{4}}\left( 4\left( \Delta x\right)
^{2}\sum_{j=1}^{J}\left( \left( u_{j}^{\prime }\right) ^{2}+\left(
u_{j-1}^{\prime }\right) ^{2}\right) \Delta x-12\left( \Delta x\right)
^{2}u_{J}u_{J}^{\prime }\right),
\end{array}
\end{eqnarray*}
which gives 
\begin{eqnarray*}
\begin{array}{l}
\displaystyle  a\left( u_{h},u_{h}\right) =\frac{k^{2}}{\left( \Delta x\right) ^{4}}\left(
3\sum_{j=1}^{J}\left( 2\left( u_{j}-u_{j-1}\right) -\Delta x\left(
u_{j}^{\prime }+u_{j-1}^{\prime }\right) \right) ^{2}\Delta x  \right.\\
\displaystyle  \left. +\left( \Delta
x^{2}\right) \sum_{j=1}^{J}\left( u_{j}^{\prime }-u_{j-1}^{\prime }\right)
^{2}\Delta x\right) .
\end{array}
\end{eqnarray*}
Hence
\begin{eqnarray*}
\begin{array}{l}
\displaystyle  a\left( u_{h},u_{h}\right)  \leq \frac{k^{2}}{\left( \Delta x\right) ^{4}}%
\left( 24\sum_{j=1}^{J}\left( u_{j}-u_{j-1}\right) ^{2}\Delta x  
\displaystyle   +6\left(
\Delta x^{2}\right) \sum_{j=1}^{J}\left( u_{j}^{\prime }+u_{j-1}^{\prime
}\right) ^{2}\Delta x  \right. \\
\displaystyle \left. +\left( \Delta x^{2}\right) \sum_{j=1}^{J}\left(
u_{j}^{\prime }-u_{j-1}^{\prime }\right) ^{2}\Delta x\right)  \\
\displaystyle  \leq \frac{24 k^{2}}{\left( \Delta x\right) ^{4}}\left(
\sum_{j=1}^{J}\left( u_{j}-u_{j-1}\right) ^{2}\Delta x+\left( \Delta
x^{2}\right) \sum_{j=1}^{J}\left( \left( u_{j}^{\prime }\right) ^{2}+\left(
u_{j-1}^{\prime }\right) ^{2}\right) \Delta x\right) 
\end{array}
\end{eqnarray*}
Using (\ref{minoration u_h}), we deduce
\begin{eqnarray*}
a\left( u_{h},u_{h}\right) \leq \frac{24 k^{2}}{\left( \Delta x\right) ^{4}}%
\frac{420 \times 19^2}{37}\left( u_{h},u_{h}\right) ,\qquad \forall u_{h}\in V_{h}
\end{eqnarray*}
and thus
\begin{equation}
\kappa \left( h\right) \le \frac{24\times 420 \times 19^2}{37}\frac{k^{2}}{%
\left( \Delta x\right) ^{4}} = \frac{2^5 \times 5 \times 19^2 \times 21}{37} \frac{k^{2}}{%
\left( \Delta x\right) ^{4}}.  \label{kappa(h)}
\end{equation}

\begin{remark}
This is certainly not an optimal upper bound for $\kappa \left( h\right) $.
\end{remark}

%
% BibTeX users please use
% \bibliographystyle{}
% \bibliography{}

\end{document}